\documentclass[a4paper]{article}
\usepackage{natbib}

  \usepackage{graphics,graphicx}
   \usepackage{amsmath}
   \usepackage{amsfonts}
   \usepackage{float,txfonts,amssymb,bm}
 \usepackage{colortbl}

\newcommand{\boldkappa}{\ensuremath{\bm{\kappa}}}

\newcommand{\boldtau}{\ensuremath{\bm{\tau}}}
\newcommand{\bh}{\ensuremath{\mathbf{h}}}
\newcommand{\br}{\ensuremath{\mathbf{r}}}

\newcommand{\tauprod}{\ensuremath{\prod_{j=1}^{N-2} \tau_j^{\rho_j}     }}

\newcommand{\taugeom}{\ensuremath{\tilde{\tau}}}
\newcommand{\tgeom}{\ensuremath{\tilde{t}}}
\newcommand{\Rho}{P}

\bmdefine\bone{1}    \bmdefine\bX{X}  \bmdefine\bt{t}
\bmdefine\btau{\tau} \bmdefine\bx{x}  \bmdefine\balpha{\alpha}
\bmdefine\bbeta{\beta}  \bmdefine\bxi{\xi}
\bmdefine\btheta{\theta}

\begin{document}

\title{Probability-Matching Predictors for Extreme
Extremes}
\author{Allan McRobie \\
Cambridge University Engineering Department\\
Trumpington St, Cambridge, CB2 1PZ, UK \\
fam20@cam.ac.uk}

\maketitle

\begin{abstract}
A location- and scale-invariant predictor is constructed which exhibits
good probability matching for extreme predictions outside the span of data drawn from
a variety of (stationary) general distributions. It is constructed via the
three-parameter $\{ \mu, \sigma, \xi \}$
Generalized Pareto Distribution (GPD). The predictor is designed to provide matching
probability exactly for the GPD in both the extreme heavy-tailed limit $\xi \rightarrow \infty$
and the extreme bounded-tail limit $\xi \rightarrow -\infty$,
whilst giving a good approximation to probability matching at all intermediate values of the tail
parameter $\xi$.
The predictor is valid even for small sample sizes $N$, even as small as $ N = 3$.

The main purpose of this paper
is to present the somewhat lengthy derivations which draw heavily on the theory of hypergeometric
functions, particularly the Lauricella functions. Whilst the construction is inspired by the Bayesian
approach to the prediction problem, it considers the case of vague prior information about both parameters
and model, and all derivations are undertaken using sampling theory.
\end{abstract}

\section{Introduction}
This paper presents a novel approach to extrapolation beyond the span of historical data, one of the basic
problems of inference. The notion of say the ``1 in 10,000 year event'' exists
in common parlance, even though the
philosophical interpretation of the phrase differs in detail between Bayesian and non-Bayesian schools. Often the analyst
has only a limited set of historical data, say of the order of 100 years, but is asked to make predictions
regarding possible future events (at say the ``1 in 10,000 year level'') which are substantially greater than historical experience. The paper concerns the issue of whether there are rational methods for tackling this philosophically-fraught problem, and the method developed here - whilst only in its preliminary stages of development - is presented as a potential candidate which may be worthy of further exploration. Although the following theory is created around properties of the Generalised Pareto Distribution (GPD) and its tail parameter $\xi$, the final predictor created is non-parametric. Its usefulness or otherwise remains to be determined, but nevertheless it has some remarkable properties.

We consider here only the univariate
case, with the data consisting of a set of $N$ discrete real-valued ``events'' $\bx$ drawn from a stationary distribution $F_{\btheta}$
with unknown parameters $\btheta$.
The approach presented focuses directly on prediction rather than estimation,
and the key concept is that of {\it probability matching}. The phrase has been given precise definitions elsewhere
(e.g. \cite{datta1, sweeting1}) and refers to cases when coverage probabilities of Bayesian and frequentist approaches coincide. 
The use of the phrase here is similar, but for the predictors in question it is not clear that a Bayesian prior exists, thus the phrase will be applied here directly to the predictor. 
Its usage here is an attempt to
formalise the loose notion that a prediction of the  ``1 in T'' level event $x_T$ should actually deliver what it appears to promise: that the probability that the next event $x_{next}$ will exceed $x_T$ is indeed $1/T$. Of course, this loose notion needs to be formalised, particularly since Bayesians will disagree with non-Bayesians about even what is meant by ``probability''.

In the frequentist approach, the analyst typically uses the data $\bx$ to construct a point {\bf estimate} $\hat{\btheta}$ of the actual unknown parameters $\btheta_0$ of the assumed model, and then
uses the tail of the family member $F_{\hat{\btheta}}$ to select the value $\hat{x}_T$, via $F_{\hat{\btheta}}(\hat{x}_T) = 1 - 1/T$. Thus $\hat{x}_T$ is an estimate of the actual (but unknown) value $x_T$
that has an exceedance probability $1/T$ for the actual (but unknown) distribution $F_{\btheta_0}$ from which the data was drawn.
Confidence intervals may then be
constructed around that estimate $\hat{x}_T$ at some chosen confidence level, and a designer or decision maker may perhaps
choose to use the upper confidence level as the basis for the ``1-in-T level'' event. Note that there is thus a mixture of notions of probability here: the designer may aim for the 1 in 10,000 event but may only have a 95\% confidence in the result.

In the Bayesian approach, the analyst combines the data with prior knowledge about the parameters to construct
the posterior distribution on the parameters, this being a representation of the analyst's updated beliefs
about the parameters of the model. Rather than simply choosing the single family member with the parameter value set to the point estimate which has the greatest posterior belief value, the Bayesian takes account of uncertainty in knowledge of the parameters by constructing the {\it predictive distribution}, this
 being the analyst's updated beliefs about the possible value of the next data element $x_{next}$. The predictive distribution
 is, loosely speaking, the integral of all probabilities of all possibilities. The Bayesian candidate for the ``1 in T'' level is that value $x_T$ above which lies a fraction $1/T$ of the analyst's beliefs about $x_{next}$.
That is, extreme value predictions at some given return level are
set at the value above which the tail integral of the predictive distribution equals
the desired exceedance probability (see \citet{colestawn}, for example). Despite the appealing rationality of the procedure, it is not without difficulties. For predictions well outside the span of the data, the tails of the predictive distribution may be strongly influenced by the prior
beliefs about model parameters, and this can be problematic when prior
knowledge - particularly regarding the tail parameters -  is vague.
Moreover, standard numerical integration techniques such as MCMC can require excessively long run
times to explore the tail regions outside the data sufficiently often for the analyst to be confident that the numerically-generated  distributions have converged sufficiently.

The problem of predicting the ``1 in T'' level event has solutions in various restricted cases. For example, the Bayesian predictor $x_T$ based on some proper prior $\Pi(\bxi)$ matches probability at the target level $T$ in the sampling sense, sampling parameters $\bxi_i$ from $\Pi(\bxi)$ and then sampling the data (and $x_{next}$) from the chosen distribution $F_{\bxi_i}$.
The two parameter ($\mu, \  \sigma$) location-scale families provide other ready examples. For example, if the data $\bx$ is sampled from some unknown member of some known family of location-scale distributions $F(y)$ with $y = (x-\mu)/\sigma$, then any statistic of the form $s_T = s_{LS} + \beta s_S$ (where $s_{LS}$ is location-scale invariant and $s_S$ is scale invariant) is a probability-matching predictor at some level $T = T(\beta)$.
Specifically, the predictor $x_T$ that arises from the Bayesian approach using the (improper)
$1/\sigma$ prior matches probability at the designed-for level $T$ (see \cite{mcrobieEEE}). (Priors using other powers of $\sigma$ give probability performance which is location- and scale-invariant, but not at the designed-for level).

Although the $1/\sigma$ predictor has the pleasing property of delivering the required $T$ level performance no matter what the parameter values actually are,
one obvious short-coming that limits its usefulness
is the requirement for {\bf complete prior knowledge}
of the model. That is, the functional form of the two-parameter family
of distributions $F$ must be known {\it a priori}. The compass of
the procedure would thus be somewhat expanded if it could be
extended to cover three-parameter location-scale-shape families of
the form $F(y)$ with $y = \left( (x-\mu)/\sigma \right)^{-1/\xi}$.
Since both the Generalized Pareto and Generalised Extreme Value
distributions (GPD, GEV) can be expressed in this form, the
possibility may then exist that the machinery of Extreme Value Theory
(e.g. \citet{Embrechts}) could also be invoked in order to apply the
predictor to data sets drawn from more general distributions. Loosely speaking,  since
 - under suitable conditions - the upper order statistics of samples drawn from
more general distributions have the GPD as their limiting distribution, then
a probability-matching predictor for the GPD may have wider application.

The question thus arises as to whether probability-matching
predictors can be constructed in the general three-parameter ($\mu$, $\sigma$, $\xi$) case.
If so, there are the further questions as to
whether there is a corresponding prior, and what form that prior might be.

This paper endeavours to construct a probability-matching predictor for the
three-parameter GPD.
By taking a sampling - rather than a Bayesian - approach,
predictors are constructed such that probability matching is exact
in the both the extreme heavy-tailed ($\xi \rightarrow \infty$) and extreme
bounded-tail ($\xi \rightarrow -\infty$)) limits.
At finite values of $\xi$, probability is only approximately matched, but the degree of approximation is very good.
Moreover the predictor applies to samples as small as $N=3$, and works remarkably well for predictions
which lie far outside the span of the data.
Finally, when applied to small data sets sampled from non-GPD distributions,
out-of-sample predictions match probability to a remarkably close approximation.

\section{Sampling Distributions of the Normalised Data}
Suppose $N$ data points $\bx = \lbrace x_1, \ldots, x_N
\rbrace $ are sampled from a Generalised Pareto Distribution (GPD) with
distribution function
\begin{equation}
F(x) = 1 - { \left(  1 + \xi \frac{(x- \mu)}{\sigma } \right) }^{-1
/ \xi}
\end{equation}
with unknown parameters $(\mu, \sigma, \xi )$. Let the ordered data
be $\bX = \mathrm{sort}(\bx)$, such that $X_1 \leq X_2 \leq \ldots
\leq X_N$. (Note that the indexing of the ordered data
is from the {\it lowest} order statistic, and that this is in the
opposite direction to that adopted in \cite{McRobieGPD} and \cite{McRobieGEV}).

The aim is to construct a predictor $x_T ( \bX ) $ such that, for
any chosen return level $T$, there is a probability $1/T$ that the
next data point $x_{next}$ will exceed $x_T$.  That is, we desire
\begin{equation}
P(x_{next} > x_T ) = \frac{1}{T}
\end{equation}
irrespective of the values of the parameters $(\mu, \sigma, \xi )$
of the distribution from which the data was sampled.

We first normalise the data to lie within the unit interval via the
statistics $\bt = \lbrace t_1, t_2, \ldots, t_{N-2} \rbrace $ with
\begin{equation}
t_j \equiv \frac{X_{j+1} - X_1}{X_N - X_1} \qquad j=1,\ldots (N-2)
\end{equation}
Clearly, $0 \leq t_1 \leq t_2 \leq \ldots \leq t_{N-2} \leq 1$.

The normalised data is location- and scale-independent, in that
$\bt(\bX) = \bt (a \bX + b \bone)$ for any $a>0$ and any $b$. The
normalisation is simply a linear mapping of the data onto the
interval $[ 0,1 ]$, the data minimum $ X_1$ mapping to zero and the
data maximum $X_N$ mapping to 1.

We are interested in extrapolating to possible large future extremes
which lie outside the span of previous data. The next
data point $x_{next}$ might not exceed the data maximum $X_N$, but we shall be most
interested in those cases when it does. We shall thus denote the next data point $x_{next}$ as $x_{N+1}$. We will likewise focus on
constructing predictions $x_T$ in that region beyond the data maximum.

The next data point $x_{next} = x_{N+1}$  and the prediction $ x_T$ may be
normalised via
\begin{equation}
s  \equiv  \frac{x_{N+1} - X_1}{X_N - X_1} \ \ \ \ \mathrm{and} \ \
\ \ s_T (\bt)  \equiv  \frac{x_T(\bX) - X_1}{X_N - X_1}
\end{equation}

Given the parameters $\lbrace \mu, \sigma, \xi \rbrace
$, the sampling density $p(\mathbf{t} | \mu, \sigma, \xi )$ of the
statistics $\mathbf{t}$ may, via the elementary integrations of
Appendix 1, be expressed in terms of a Lauricella function $F_D$.
Since the statistics $\bt$ are location- and scale-independent, the
sampling density retains a parameter dependence only through the
shape parameter $\xi$.

Writing $N_1 = N-1$ for brevity, then in the region of the
heavy-tailed ($\xi > 0$) GPDs, writing $\alpha \equiv \xi^{-1}
>0 $, we obtain from Appendix 1 that
\begin{eqnarray}
p(\bt | \xi, (\xi>0) ) & = & N_1! \ \alpha^{N_1} \ \ \Gamma \left[
\begin{array}{c} N_1, \ N_1 \alpha \\ N_1(1+\alpha) \end{array} \right]
  \ F_D^{(N-2)}\left( N_1, \bm{1}+\bm{\alpha};
N_1(1+\alpha); \boldtau \right) \label{palpha}
\end{eqnarray}
where  $\bone$ is a vector of $(N-2)$ ones, $\bm{\alpha} \equiv
\alpha \bone$ and $\btau \equiv \bone-\bt$.

For the bounded-tail ($\xi < 0$) GPDs, writing $\beta \equiv
-\xi^{-1} > 0 $, we similarly obtain
\begin{eqnarray}
p(\bt | \xi, (\xi<0) ) & = & N_1! \ \beta^{N_1} \ \ \Gamma \left[
\begin{array}{c} N_1, \ N_1 \beta \\ N_1+\beta \end{array} \right]
  \ F_D^{(N-2)}\left( N_1, \bm{1}-\bm{\beta},
N_1+\beta; \bt \right) \label{pbeta}
\end{eqnarray}
with $\bm{\beta} \equiv \beta \bm{1}$.

\section{The Problem Statement}
Since all dependence on the location and scale parameters has been
removed by the normalisation, the only parameter of remaining
interest is the tail parameter $\xi$.

Given the parameters $\xi$, the probability that the next
(normalised) event $s$ will exceed some given function $s_T(\bt)$ of
any data $\bt$ is given by
\begin{equation}
P(s > s_T(\bt)| \xi) = \int \bm{G}(s_T(\bt)|\xi) \ p(\bt | \xi) \ d\bt
\label{eqnpp}
\end{equation}
where $d\bt = dt_1 dt_2 \ldots dt_{N-2}$ and the integral is over
all admissible normalised data $\bt$, and the integrand is defined
via
\begin{equation}
G(s_T(\bt)|\xi)p(\bt|\xi) \equiv \int_{s_T(\bt)}^\infty \
p(s,\bt|\xi) \ ds
\end{equation}

In the heavy-tailed region ($\alpha \equiv 1/\xi > 0$), for the case of
interest where $s_T > 1$, the
 elementary integrations of Appendix 2 reveal the integrand of equation \ref{eqnpp} to be
\begin{eqnarray}
\bm{G}(s_T(\bt)|\xi) \ p(\bt | \xi) & = & \frac{N!}{N+1} \
\alpha^{N_1} \ \ \Gamma \left[
\begin{array}{c} N_1, \ N \alpha \\ N_1+N\alpha \end{array} \right]
\frac{t_1^{N_1 \alpha -1 }}{\left( t_2 \ldots
t_{N-2}\right)^{1+\alpha} s_T^\alpha}  \ldots \nonumber \\
 &  & \hspace*{10mm} \times \ F_D^{(N-1)}\left( N\alpha,  \bm{1}+ \bm{\alpha}, \alpha ; N_1 + N\alpha ; \boldkappa  \right)
 \label{eqnpp1}
\end{eqnarray}
where the $j$-th element of $\boldkappa$  is $\kappa_j = 1 - t_1 /
t_j $ for $j = 2, \ldots, N$, the $t_j$ notation being extended here
to include the two further points $t_{N-1} = 1$ and $t_N = s_T$.

For the bounded-tail region  ($\beta \equiv -1/\xi > 0$),  the
corresponding integrand is
\begin{eqnarray}
\bm{G}(s_T(\bt)|\xi) \ p(\bt | \xi) & = & \frac{N!}{N+1} \
\beta^{N_1} \ \ \Gamma \left[
\begin{array}{c} N_1, \ 1+\beta \\ N+\beta \end{array} \right]
 \ \frac{1}{ s_T^{N_1}}  \ldots \nonumber \\
 &  & \hspace*{7mm} \times \ F_D^{(N-1)}
 \left( N_1,  \bm{1}- \bm{\beta}, 1-\beta ; N+\beta ; \frac{\bm{t}}{s_T}, \frac{1}{s_T}  \right)
\label{eqnpp2}
\end{eqnarray}

The problem statement is thus: for any $T > N+1$, find a function
$s_T(\bt)$ over the normalised data space which, when substituted
into Eqns \ref{eqnpp}-\ref{eqnpp2},  is such that $P(s> s_T(\bt)|\xi) =
1/T$ for any $\xi$.

\pagebreak

\section{The Heavy-Tailed Limit $\xi \rightarrow \infty$}

We consider first those predictors which guarantee to match probability
at the $1/T$ level for those distributions in the
limit of extremely heavy tails $ \xi \rightarrow \infty$.
For large $\xi$ (i.e.\ small $\alpha$), the
Lauricella function in Eqn.~(\ref{eqnpp1}) approaches unity due to
the argument $N \alpha$ in its first slot. The conditional
exceedance integrand of Eqn.~(\ref{eqnpp1}) thus simplifies to
\begin{equation}
 \bm{G}(s_T(\bt)|\xi) p (\bt|\xi) \approx \frac{N!}{N+1} \ \alpha^{N_1} \
\Gamma \left[
\begin{array}{c} N_1, \ N\alpha \\ N_1+N\alpha \end{array} \right]
 \
\frac{t_1^{N_1 \alpha -1 }}{\left( t_2 \ldots
t_{N-2}\right)^{1+\alpha} s_T^\alpha} \label{eqnpp3}
\end{equation}

Since any function of the statistics $\bt$ is a predictor, there is
an almost limitless variety to the possible functional forms that
our predictor $s_T(\bt)$ might take.

To progress, we consider predictors taking a power law form:
\begin{equation}
s_T = \prod_{j=1}^{N-2} \ {t_j}^{-\lambda_j} 
\label{powerlaw}
\end{equation}

Integrating (\ref{eqnpp3}) over the domain $0 \leq t_1 \leq \ldots
\leq t_{N-2} \leq 1$ and setting the result equal to the desired
exceedance probability $1/T$ leads to the constraint

\begin{equation}
\prod_{j=1}^{N-2} \left( 1 + \frac{\gamma_j}{ N-j} \right)
 = \frac{T}{N+1} \ \ \mathrm{where} \ \ \gamma_j = \sum_{k=1}^j \lambda_k \label{constraint}
\end{equation}

For $N=3$, this constraint requires
\begin{equation}
\lambda_1 = \gamma_1 = 2 \left( \frac{T}{4}- 1 \right)
\label{predN3}
\end{equation}

For larger samples and for a given $T$, the constraint equation
(\ref{constraint}) defines an $(N-3)$-dimensional space of possible
exponents for our power-law predictor. One obvious solution sets
each term in the left-hand product of equation~(\ref{constraint}) equal to the same value,
$(T/(N+1))^{1/(N-2)}$, giving
\begin{equation}
\gamma_j = (N-j)\left[ \left( \frac{T}{N+1}\right)^{\frac{1}{N-2}}
-1 \right] \label{weib1}
\end{equation}

An alternative choice is one that gives predictions that are in some
sense small, and the expected value of $\log(s_T)$ can be minimised
by choosing
\begin{equation}
\gamma_j = (N-j-1)K -(N-j) \ \ \mathrm{with} \ \ K =
\left[\frac{N-1}{N+1} \ T  \right]^{\frac{1}{N-2}} \label{weib2}
\end{equation}

For simplicity, we might instead choose to make all the $\lambda_j$
the same ($= \lambda$). For any sample size $N$ and any desired
return level $T$, the exponents of the power-law predictor can be
obtained numerically (e.g.\ by a simple bisection method) to
determine that value of $\lambda$ which satisfies the constraint
equation. Values of $\lambda$ so determined are given in Table
\ref{lambdas}.

\begin{table}[h!]
\caption{Exponents $\lambda$ satisfying the heavy-tailed constraint
equation \ref{constraint} \label{lambdas}} \centering
\begin{tabular}{|r||r|c|c|c|}
  \hline
 \rule[-0.2cm]{0cm}{0.8cm}  $\frac{T}{N+1}$ & N=3 & N=7 & N=15 & N=31 \\ \hline \hline
    2 &    2 & 0.1507 & 0.0354 & 0.0113 \\
    4 &    6 & 0.3363 & 0.0750 & 0.0233 \\
    8 &   14 & 0.5615 & 0.1188 & 0.0361 \\
   16 &   30 & 0.8317 & 0.1672 & 0.0498 \\
   32 &   62 & 1.1530 & 0.2202 & 0.0642 \\
   64 &  126 & 1.5325 & 0.2782 & 0.0795 \\
  128 &  254 & 1.9785 & 0.3412 & 0.0956 \\
  256 &  510 & 2.5003 & 0.4097 & 0.1125 \\
  512 & 1022 & 3.1086 & 0.4838 & 0.1303 \\
 1024 & 2046 & 3.8157 & 0.5638 & 0.1489 \\
 2048 & 4094 & 4.6359 & 0.6501 & 0.1684 \\
 4096 & 8190 & 5.5854 & 0.7430 & 0.1888 \\
  \hline
\end{tabular}
\end{table}

\begin{figure}[h!] \centering
  \includegraphics[width=120mm,keepaspectratio]{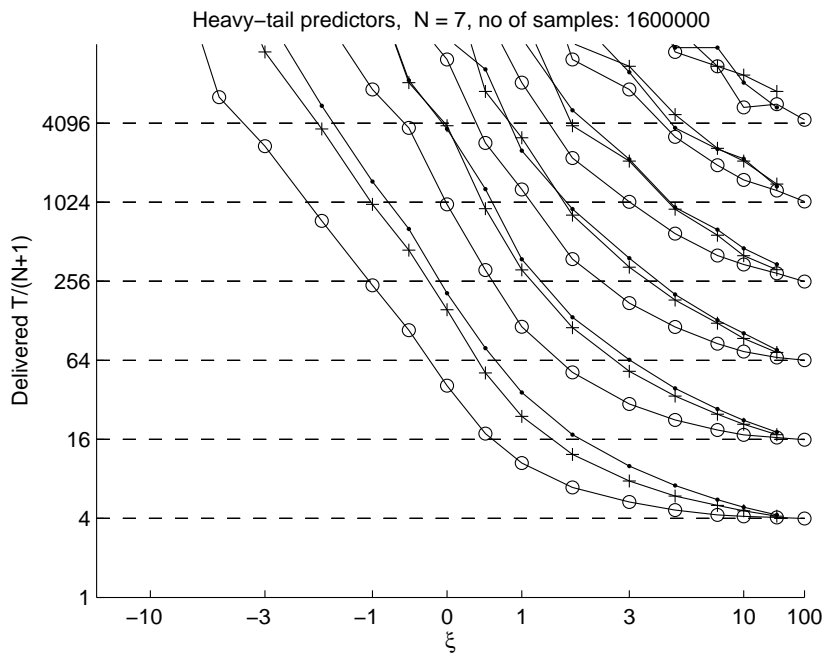}
  \caption{ The probability performance of the heavy-tailed predictors.
The figure shows the return levels delivered by the predictors of
Eqns. \ref{weib1} (.), \ref{weib2}(+) and with all $\lambda_j =
\lambda$ (o) for samples of size $N = 7$.
In each case, predictors were asked to predict at return levels of $T/(N+1) = 4, 16, 64,256, 1024 \ \text{and} \ 4096$.
All three predictors
satisfy the heavy-tail constraint equation \ref{constraint}
and the convergent lines at the right-hand side of the figure support the proposition that each
matches probability at the target level in the heavy-tailed limit.
  }\label{threeweib}
\end{figure}

The performance of the above predictors when played against samples
drawn from distributions with various shape parameters $\xi$ are
shown in Figure \ref{threeweib}. All give the desired performance in
the extremely heavy-tailed limit (at the right of the diagram), but over-predict elsewhere.

The predictors in this heavy-tailed limit we shall denote by $u_\alpha = s_T -1$,
this being the (scaled) excess of the prediction above the data maximum.

\vspace*{30mm}

\section{The Bounded-Tail Limit $\xi \rightarrow -\infty$}
In this section a constraint equation is constructed for the
exponents of a power law predictor which will deliver exact
probability matching in the extreme limit of large
negative $\xi$.  Such a constraint is somewhat more difficult to
find. The derivation in Appendix 3 shows that a
predictor can be constructed of the form $ s_T = 1+u_\beta$ where
\begin{equation}
u_\beta = \tau_{N-2} \ \frac{\btau^\Rho}{1-\btau^\Rho} \ \
\mathrm{where} \ \btau^\Rho \equiv \prod_{j=1}^{N-2} \tau_j^{\rho_j}
\label{ufpred}
\end{equation}
 and the exponents $\rho_j$ of each $\tau_j = 1-t_j$ must satisfy the
constraint
\begin{equation}
\prod_{k=1}^{N-2} \left[ \frac{2\eta_k}{k+2}+1 \right] =
\frac{1}{1-\frac{N+1}{T}} \ \mathrm{where} \ \eta_k = \sum_{j =
N-1-k}^{N-2} \rho_j \label{boundedconstraint}
\end{equation}

For $N=3$ there is a single solution
\begin{equation}
\rho_1 = \eta_1 = \frac{3}{2 \left( \frac{T}{4}-1 \right) } =
\frac{3}{\lambda_1}
\end{equation}
which is reciprocally related to the exponent $\lambda_1$ of the
corresponding extreme heavy-tailed case.

Again, for larger sample sizes and given $T$, the constraint defines
an $(N-3)$-dimensional manifold of possible exponents, and various
ad hoc schemes can be readily devised that satisfy the constraint.
For further progress, we consider only the scheme which sets all
exponents $\rho_j$ equal to the same value $\rho$ for some given
$T$. Since the constraint equation is a polynomial of degree $N-2$
in $\rho$, this can again be solved numerically for $\rho$ at any
$T$, for example by using a simple bisection method. Values of
$\rho$ for $N=3$, 7, 15 and 31 are given in Table~\ref{myrhos}.

\begin{table}[h!]
\caption{Exponents $\rho$ satisfying the extreme bounded-tail constraint equation
\ref{boundedconstraint} \label{myrhos}} \centering
\begin{tabular}{|r||r|c|c|c|}
  \hline
 \rule[-0.2cm]{0cm}{0.8cm}  $\frac{T}{N+1}$ & N=3 & N=7 & N=15 & N=31 \\ \hline \hline
    2 &    1.5                  & 0.1326 & 0.0381 & 0.0147 \\
    4 &    0.5                  & 0.0527 & 0.0155 & 6.0377$\times 10^{-3}$ \\
    8 &   0.2066                & 0.0241 & 7.1688$\times 10^{-3}$& 2.7948$\times 10^{-3}$ \\
   16 &   0.1                   & 11.544$\times 10^{-3}$ & 3.4548$\times 10^{-3}$ & 1.3491$\times 10^{-3}$ \\
   32 &   0.0484                & 5.6596$\times 10^{-3}$ & 1.6975$\times 10^{-3}$ & 0.6633$\times 10^{-3}$ \\
   64 &  0.0238                 & 2.8022$\times 10^{-3}$ & 0.8414$\times 10^{-3}$ & 0.3271$\times 10^{-3}$ \\
  128 &  11.809$\times 10^{-3}$ & 1.3945$\times 10^{-3}$ & 0.4189$\times 10^{-3}$ & 0.1638$\times 10^{-3}$ \\
  256 &  5.8796$\times 10^{-3}$ & 0.6953$\times 10^{-3}$ & 0.2089$\times 10^{-3}$ & 0.0817$\times 10^{-3}$ \\
  512 & 2.9337$\times 10^{-3}$  & 0.3457$\times 10^{-3}$& 0.1044$\times 10^{-3}$ & 0.0408$\times 10^{-3}$ \\
 1024 & 1.4644$\times 10^{-3}$ & 0.1733$\times 10^{-3}$ & 0.0521$\times 10^{-3}$ & 0.0204$\times 10^{-3}$ \\
 2048 & 0.7186$\times 10^{-3}$ & 0.0866$\times 10^{-3}$ & 0.0261$\times 10^{-3}$ & 0.0102$\times 10^{-3}$ \\
 4096 & 0.3662$\times 10^{-3}$ & 0.0432$\times 10^{-3}$ & 0.0130$\times 10^{-3}$ & 0.0050$\times 10^{-3}$ \\
  \hline
\end{tabular}
\end{table}

\begin{figure}[h!] \centering
  \includegraphics[width=120mm,keepaspectratio]{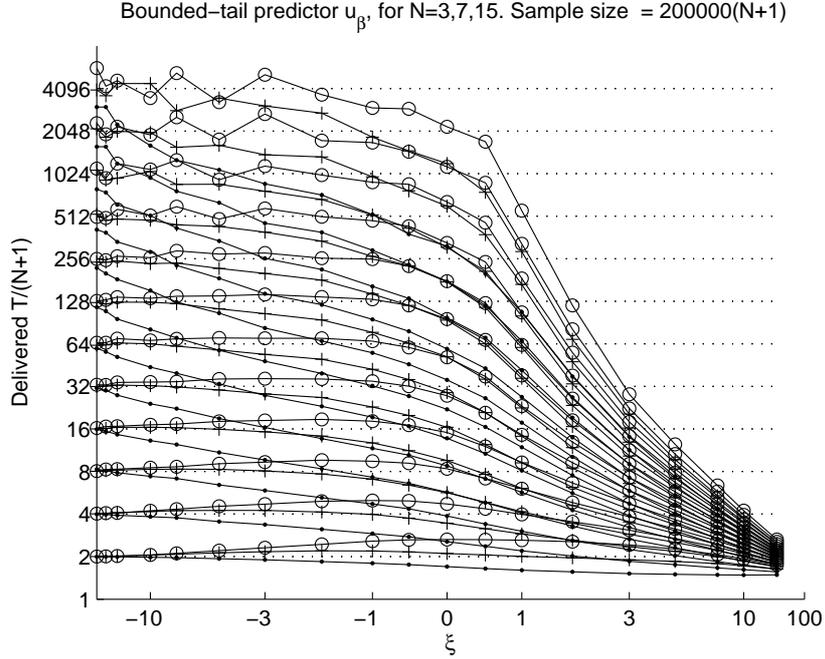}
  \caption{ The probability performance of the extreme bounded-tail predictor $u_\beta$,
  wherein all exponents
$\rho_j = \rho$ satisfy the constraint Eqn. \ref{boundedconstraint} derived
in Appendix 3. The figure gives the return levels delivered for
samples of size $N=3$(.), 7(+) and 15(o), and the convergent lines to the left of the figure
suggest that such
predictors do indeed match probability in the extreme bounded-tail limit.
}\label{threebounded}
\end{figure}

The performance of the resulting predictor is illustrated in
Fig.~\ref{threebounded} for $N=3$, 7 and 15, showing that $1/T$
exceedance is delivered in the extreme bounded-tail limit $\xi \rightarrow -\infty$, but typically
under-predicting for less extreme tail parameters.

\pagebreak

\section{Predictors at specific intermediate values of $\xi$}
We now possess predictors $u_\alpha$ and $u_\beta$ which match probability in the two extreme limits.
Individual predictors can readily be constructed at {\bf any specific
known} $\xi$ between these two limits. For example at the GPDs where
$\xi = 0$ and $\xi = -1$ (the exponential and the uniform) any
number of location- and scale-invariant predictors can be
constructed which match probability there, and a simple example is provided by the Bayesian
$1/\sigma$ predictors.

For $\xi = 0$ (exponential), the $1/\sigma$ Bayesian predictor is
\begin{equation}
s_T = (1+ \sum_{j=1}^{N-2}t_j) \ \left[ \left( \frac{NT}{N+1}
\right)^{1/(N-1)} -1 \right] \label{exppred}
\end{equation}

and for $\xi = -1$ (uniform) it is simply the constant
\begin{equation}
s_T = \left( \frac{T}{N+1} \right)^{1/(N-1)} \label{unifpred}
\end{equation}

\begin{figure}[h!] \centering
  \includegraphics[width=70mm,keepaspectratio]{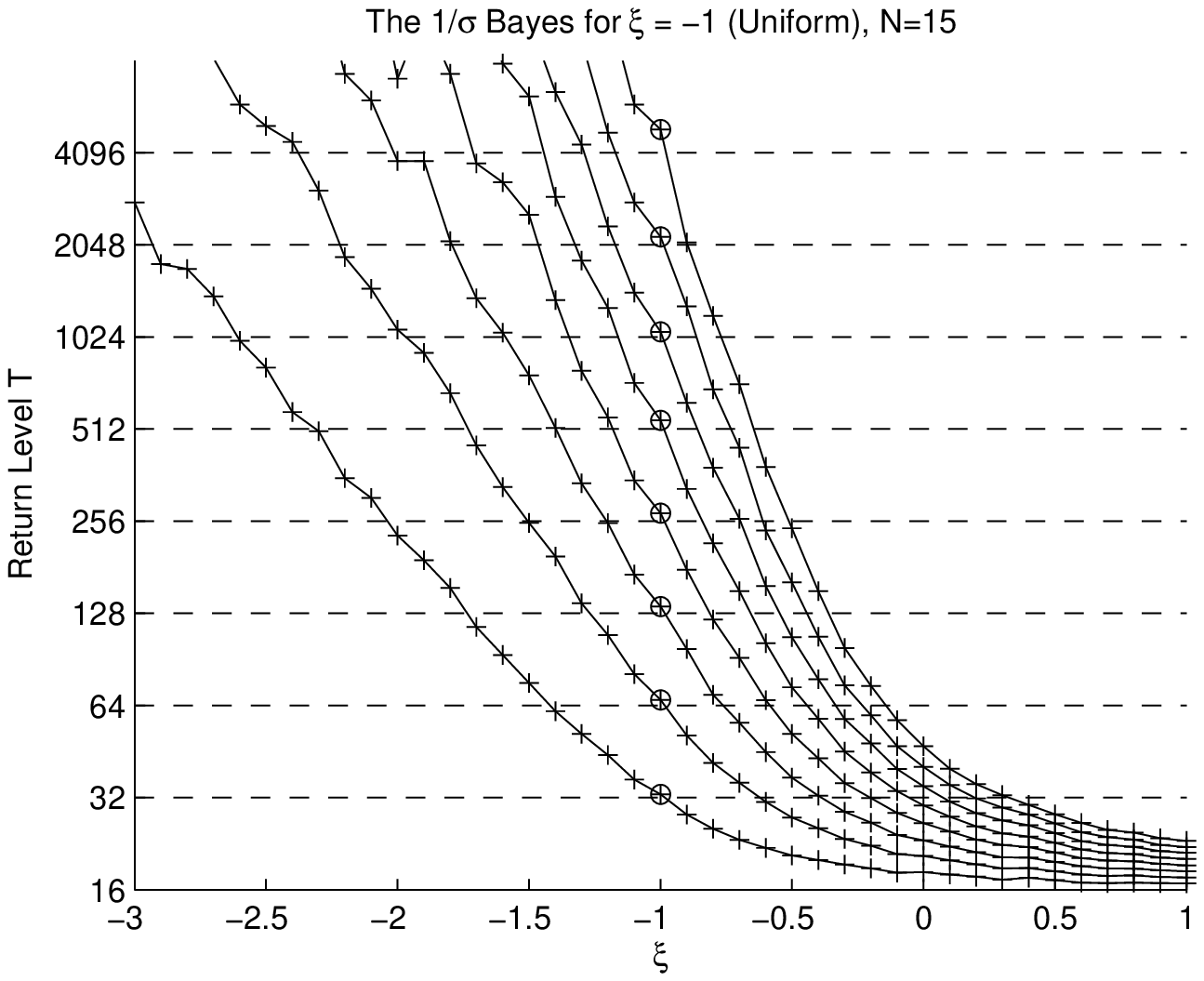}
 \includegraphics[width=70mm,keepaspectratio]{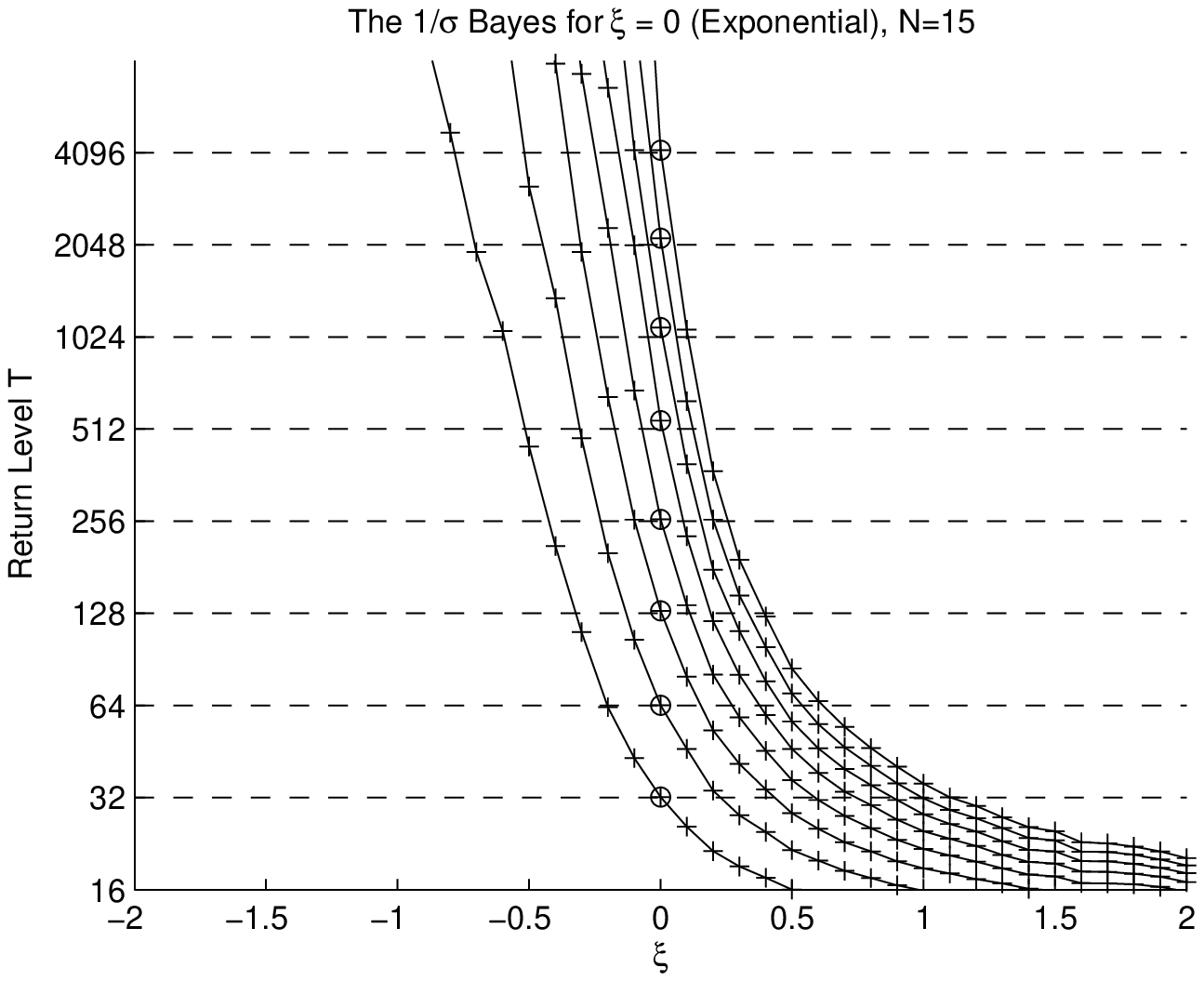}
  \caption{The numerical performance of the $1/\sigma$ Bayes
  predictors for $\xi = -1$ uniform (left) and $\xi = 0$ exponential (right).
  Both deliver the desired exceedance performance only at the respective shape
  parameters (along the central vertical in each diagram) that they have been designed for.
  }\label{xizeroone}
\end{figure}

The performance of these two predictors over a range of shape
parameters is illustrated in Fig.\ \ref{xizeroone}. It can be seen
that they do indeed deliver the required prediction performance at
the shape parameter for which they are designed, but the performance
deviates rapidly away from the desired level at nearby shape
parameters. This illustrates the danger of designers assuming that
data is say exponentially distributed and predicting accordingly,
for if the data is actually drawn from a nearby GPD the predictions
might be extremely optimistic, and the designer should not be so
surprised when the design is soon exceeded.

Note that the exponential $1/\sigma$ predictor does not predict
the data maximum ($s_T = 1$) for the return level $T = N+1$. It can even predict
values below the data maximum (e.g. for $N= 3$, $s_{T=4} = (\sqrt{3} -1)(1+t)$, which can be as low as $s_{T=4} = 0.7321$).
This
differs from the strategy that will be adopted in this paper. Here
the aim is to create more general predictors and the decision has
thus been made that these should all pass through the only presently
known (albeit trivial) example of a universal probability-matching
extreme value predictor (namely the data maximum for the $T= N+1$
prediction, corresponding to $s_{T=N+1} = 1$).

In the absence of a known non-informative prior on the shape
parameter, a number strategies could be adopted for constructing
predictors which attempt to match probability across the full range of $\xi$.
In this paper, the approach will be to interpolate between the two extreme predictors
$u_\alpha$ and $u_\beta$.

\section{Predictions across all shape parameters}
The two extreme predictors proposed thus far each define a
surface $u(\bt)$ above the $(N-2)$-dimensional simplex of all
possible normalised data $\bt$. If there is a general predictor
which matches probability at all $\xi$, then one might expect its
prediction surface $u(\bt)$ to approach $u_\alpha(\bt)$ for small $t_j$
and $u_\beta(\bt)$ for $t_j$ near unity, these being the regions where
data tends to congregate in the respective extreme limits.

A large number of interpolation schemes were considered
but for brevity, only one such scheme is presented here.

\subsection{An interpolated predictor}
First we pre-condition the two extreme predictors $u_\alpha$ and $u_\beta$ such that
they give better probability matching over wider ranges of $\xi$ than the extreme limits for which
they have been designed, and then we interpolate.

\begin{figure}[h!] \centering
  \includegraphics[width=70mm,keepaspectratio]{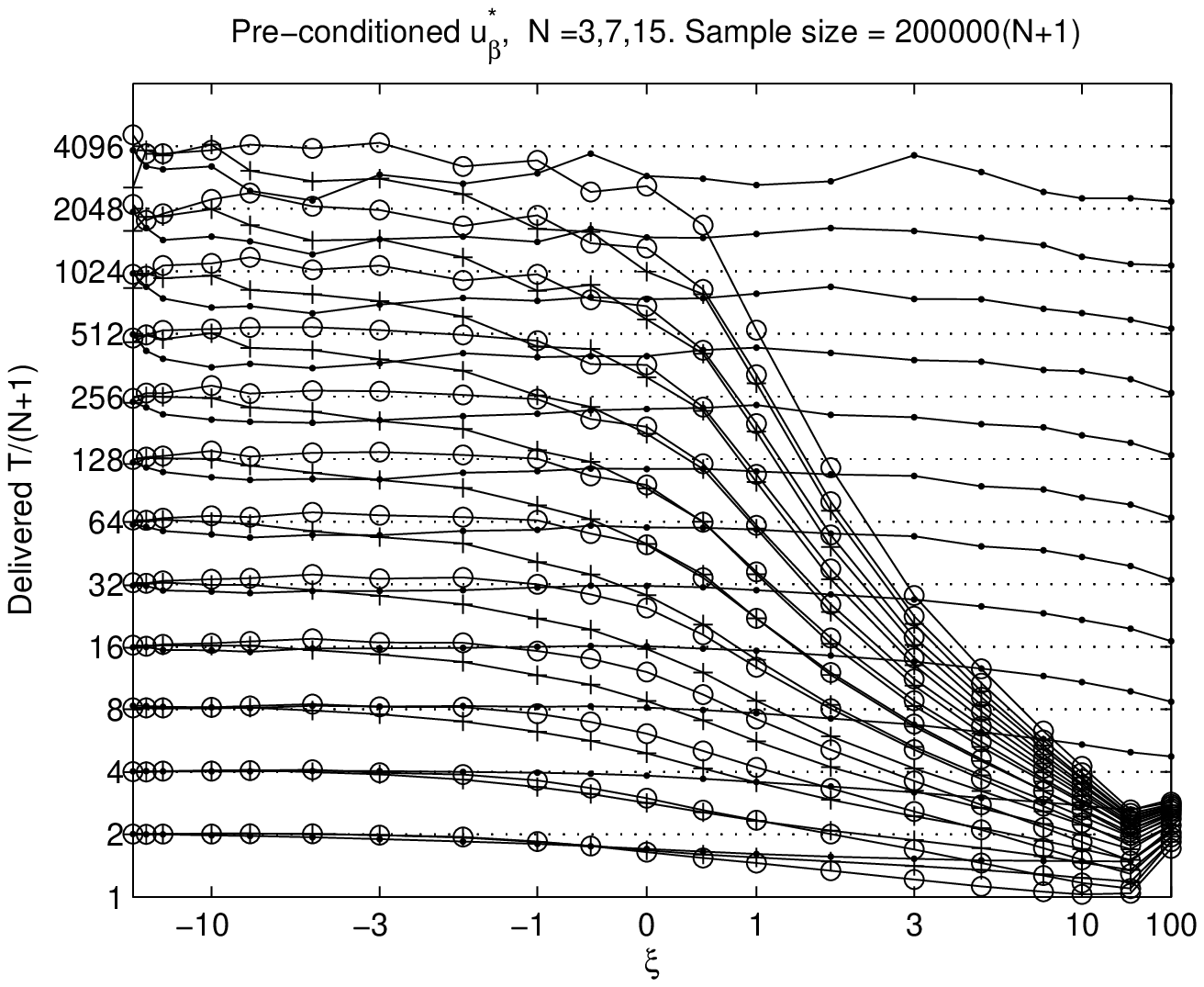}
\includegraphics[width=70mm,keepaspectratio]{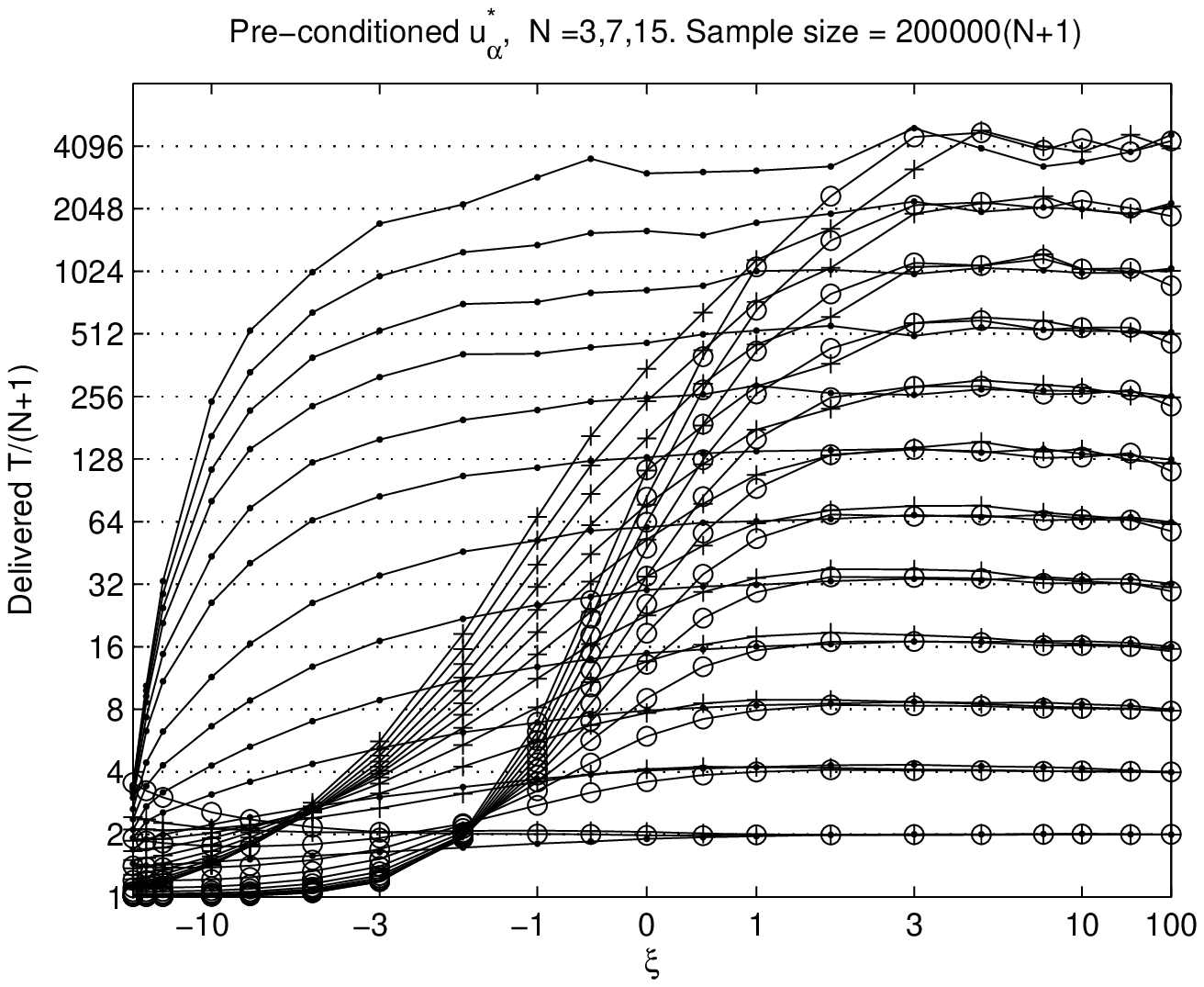}
  \caption{The probability performance of the pre-conditioned predictors
$u_\beta^*$ and $u_\alpha^*$  for samples of size $N = 3$ (.), 7 (+)  and 15 (o), showing reasonable
probability matching over large areas of the respective $\xi < 0$ or $\xi >0$ regions.
  }\label{precon}
\end{figure}

Pre-multiplying $u_\alpha$ and $u_\beta$  by
some power of the geometric means $\tilde{\tau}$, $\tilde{t}$ of
$\btau$, $\bt$ respectively, can improve their performance over
wider ranges of $\xi$ whilst matching probability in their respective limits.
That is, moderated predictors can be constructed of the form
\begin{eqnarray}
u_{\alpha}^{*} & = & \taugeom^{A} u_\alpha \\
u_{\beta}^{*} & = & \tgeom^{B} u_\beta   \label{AB}
\end{eqnarray}
where the exponents $A$ and $B$ are chosen by numerical experiments
such that $u_\alpha^*$ and $u_\beta^*$ give reasonable probability matching
over most of their respective $\xi>0$ and $\xi <0$ ranges (see Fig.~\ref{precon}).
Candidate values for $A$ and $B$ determined on the basis of such
numerical experiments are shown in Table~\ref{ABtable}.

\begin{table}[h!]
\caption{Empirically-determined exponents $A$ and $B$ of the
moderated predictors of Eqn. \ref{AB} used to construct Figs. \ref{precon} and
\ref{may096} \label{ABtable}} \centering
\begin{tabular}{|r||c|c|c|c||c|c|c|c|}
  \hline
& \multicolumn{4}{c||}{A} & \multicolumn{4}{c|}{B} \\
  \hline
 \rule[-0.2cm]{0cm}{0.8cm}  $\frac{T}{N+1}$ & N=3 & N=7 & N=15 & N=31 & N=3 & N=7 & N=15 & N=31 \\ \hline \hline
    2 &    4 &  2.2  & 2.5 &  3   & 0    & 0.45  & 0.7   & 0.75\\
    4 &    2 &  2.38 & 3.0 & 3.5  & -2   & 0.3   & 0.5   & 0.55\\
    8 &  1.5 &  2.57 & 3.5 & 4.2  & -6   & 0.2   & 0.3   & 0.4\\
   16 & 1.25 &  2.78 & 4.0 & 5.05 & -14  & 0.1   & 0.2   & 0.3\\
   32 &   1  &  3.02 & 4.5 & 6    & -30  & 0.05  & 0.1   & 0.2\\
   64 &  0.8 &  3.3  & 5.0 & 7    & -62  & 0.02  & 0.05  & 0.1\\
  128 &  0.6 &  3.6  & 5.5 & 8    & -126 & 0.01  & 0.25  & 0.05\\
  256 & 0.55 &  3.9  & 6.0 & 9    & -254 & 0.005 & 0.0125& 0.025\\
  512 & 0.5  &  4.2  & 6.5 & 10   & -510 & 0.0025& 0.0063& 0.0125\\
 1024 & 0.5  &  4.5  & 7.0 & 11   & -1022& 0.0012& 0.0031& 0.0063\\
 2048 & 0.5  &  4.8  & 7.5 & 12   & -2046& 0.0006& 0.0016& 0.0031\\
 4096 & 0.5  &  5.1  & 8.0 & 13   & -4094& 0.0003& 0.0008& 0.0016\\
  \hline
\end{tabular}
\end{table}

A combined predictor can then be constructed via simple linear
interpolation between the two moderated prediction surfaces. The
interpolation chosen is based on the ``elemental estimators'' of \cite{McRobieGPD}.
These are a family of simple location- and scale-invariant estimators
based on log-spacings of the data, and they were shown to be absolutely unbiased estimators of the
shape parameter $\xi$ of the GPD. The specific estimator $\hat{\xi}$ used here is the one which gives equal weight
to each elemental estimator.

The interpolation functions are
\begin{equation}
f_1 = \frac{e^{\hat{\xi}}}{1+e^{\hat{\xi}}} \ \ \ \mathrm{and} \
\ f_2 = 1-f_1 = \frac{e^{-\hat{\xi}}}{1+e^{-\hat{\xi}}}
\end{equation}

The resulting total predictor $u_{T}$ is
\begin{equation}
u_{T} = f_1 u_{\alpha}^{*} + f_2 u_{\beta}^{*} \label{hybrid}
\end{equation}

All elements of this predictor have some analytical justification
except for the two numerically-determined pre-conditioning exponents $A$ and $B$.
Prediction at the return level $T$ from
any data sample of size $N$ can thus be computed by the above
formula, and requires knowledge of only two numbers, $A$ and $B$,
determined by numerical experiments.

\begin{figure}[t!] \centering
\includegraphics[width=70mm,keepaspectratio]{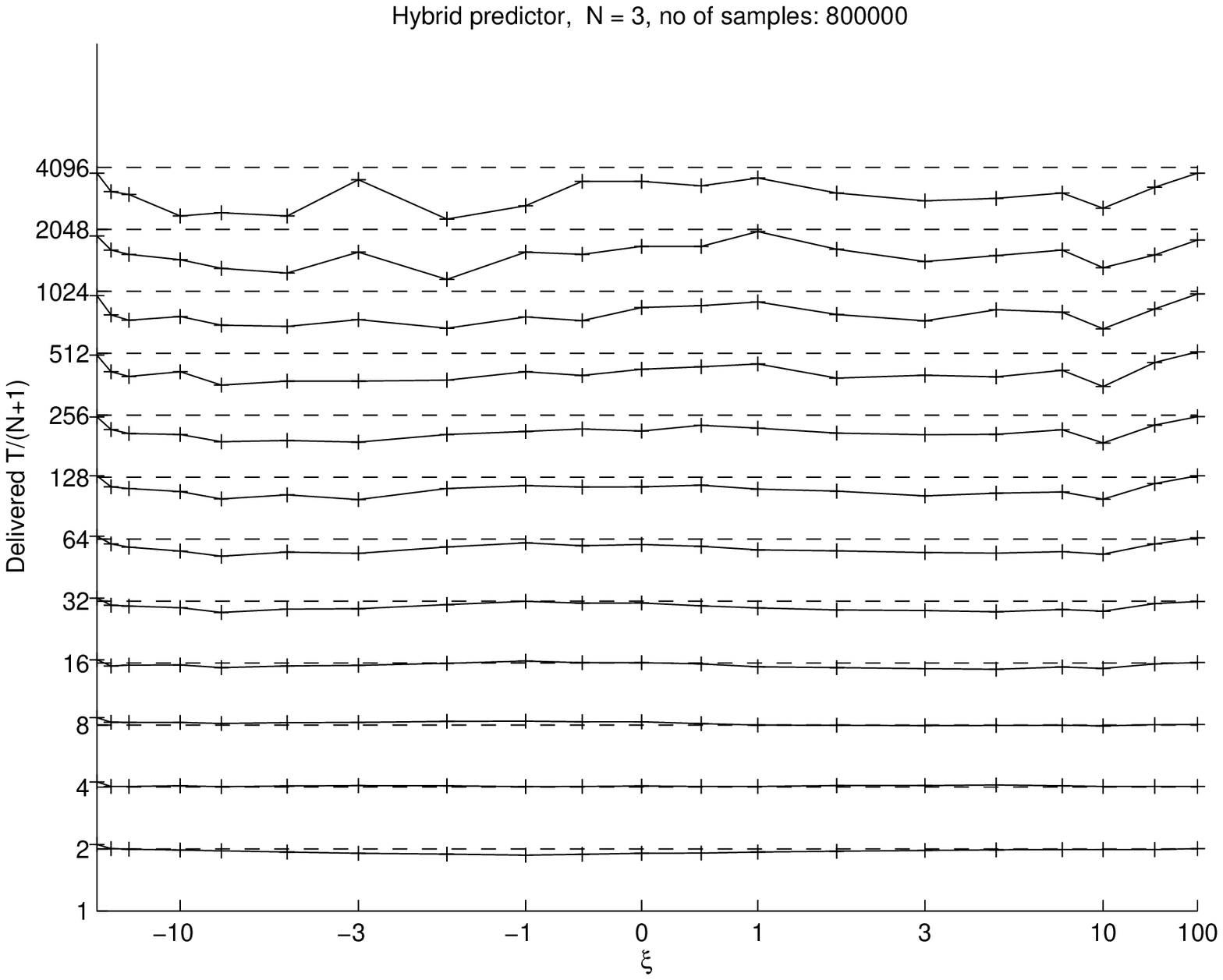}
\includegraphics[width=70mm,keepaspectratio]{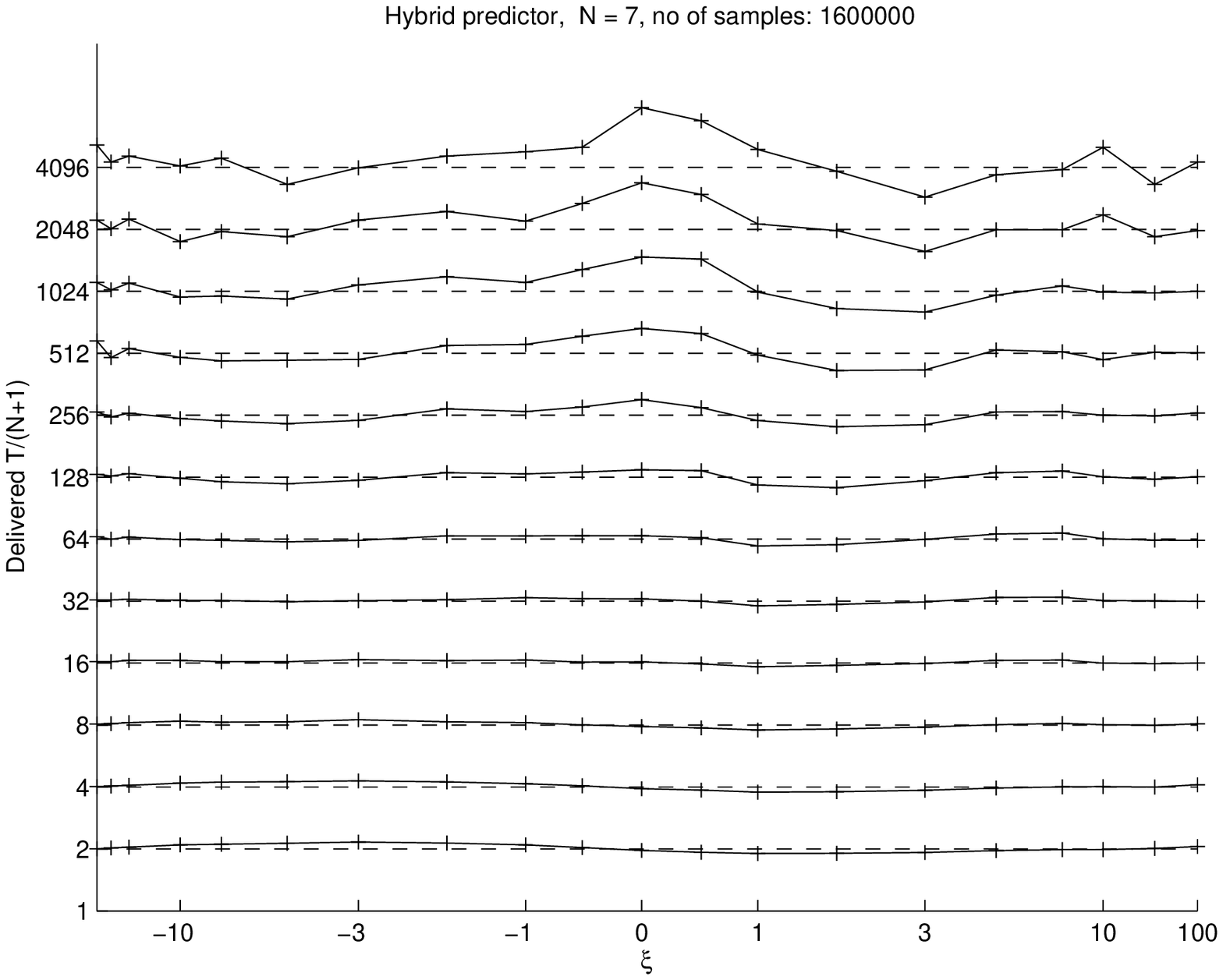}\\
\includegraphics[width=70mm,keepaspectratio]{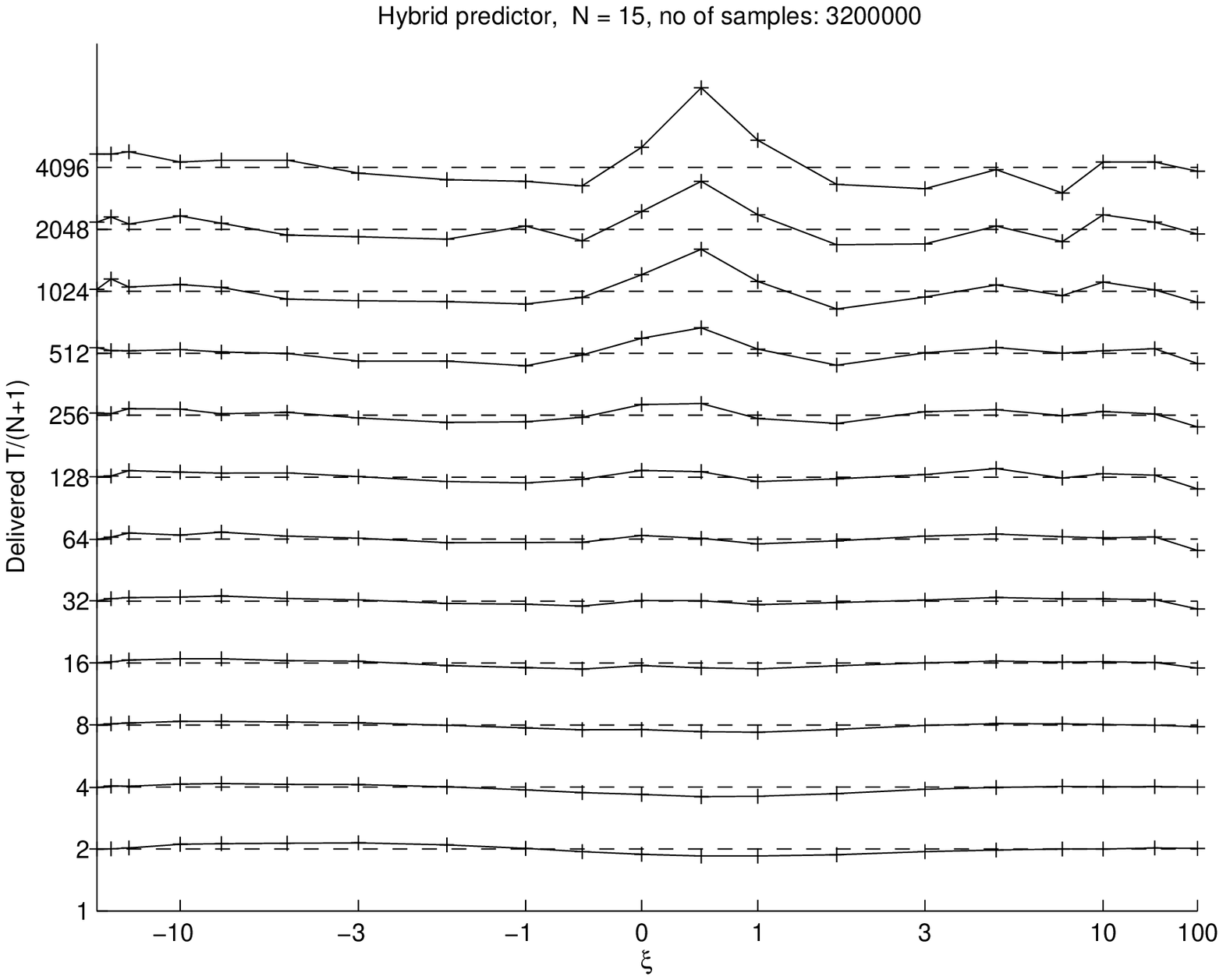}
\includegraphics[width=70mm,keepaspectratio]{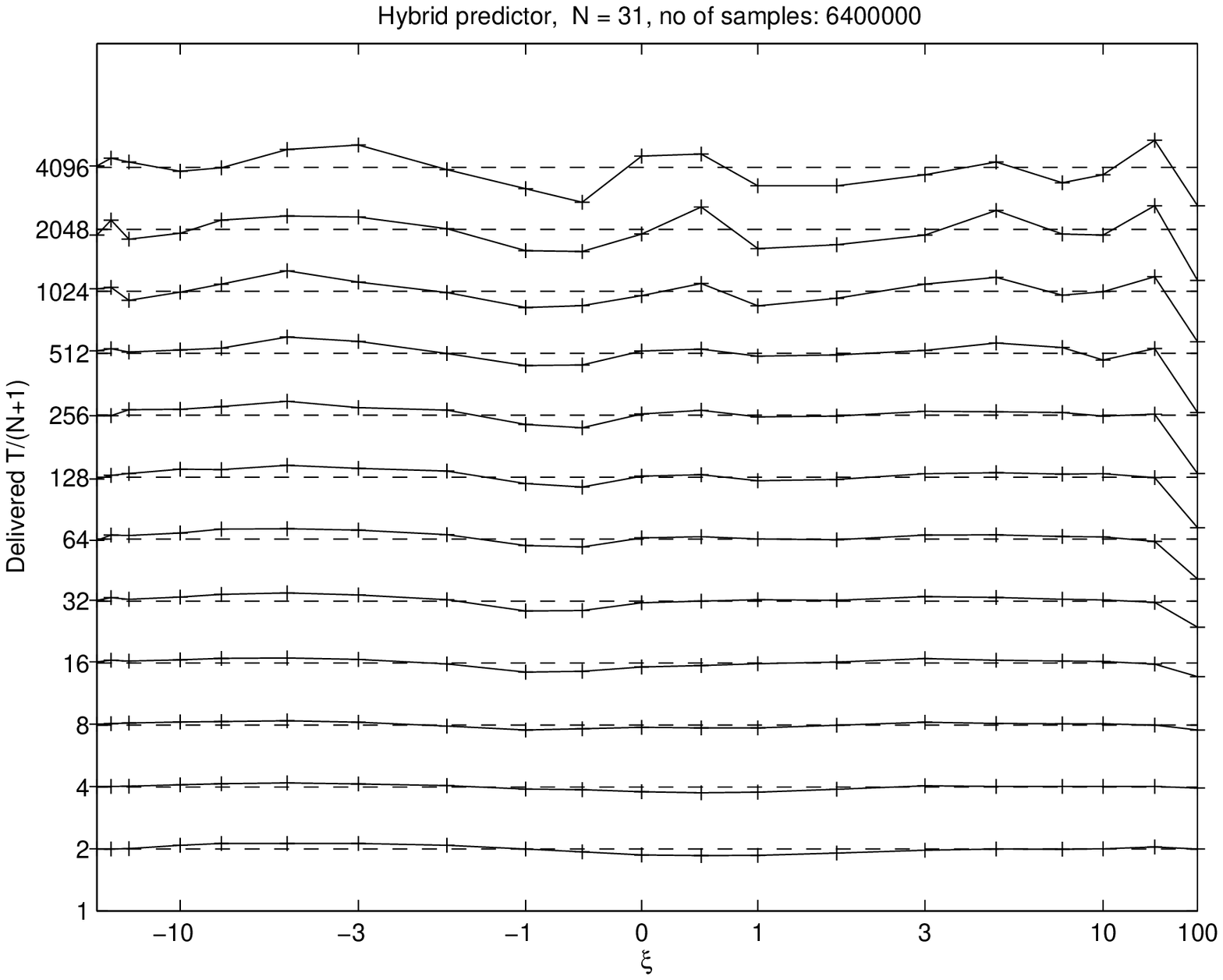}
  \caption{The probability performance of the interpolated predictor
of Eqn. \ref{hybrid} for samples of size $N = 3$, 7, 15 and 31. Good
probability matching is obtained for extrapolation factors
$T/(N+1)$ of up to 100 or more.
  }\label{may096}
\end{figure}

The performance of the interpolated predictor is shown in Fig.
\ref{may096} for samples of size $N=3$, 7, 15 and 31. It can be seen
that probability is matched to a good approximation in all cases across the full range of $\xi$,
even for extrapolations to return levels far beyond the span of the
data. The probability matching is almost exact in the lower half
of each figure, which corresponds to prediction factors $T/(N+1)$ of
up to 100 - i.e. a predictor with a return level up to $T \approx 300$
can be constructed for a sample of size 3, and up to $T \approx
3000$ for a sample of size 31. Such large extrapolations beyond the
data are often demanded in engineering: design for the 10,000 year
event often being required from around 100 years of historical data.

\begin{figure}[t!] \centering
  \includegraphics[width=100mm,keepaspectratio]{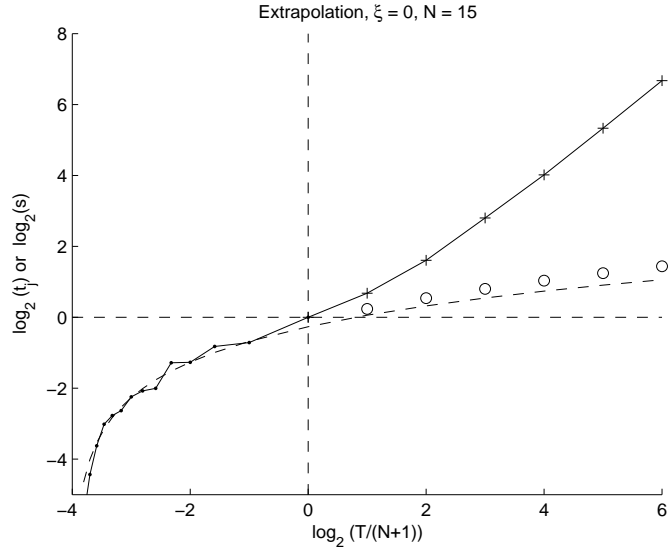}
  \caption{Typical extrapolation by the interpolated predictor. The data, shown by dots in the bottom left corner, is a random sample of size $N=15$ from an exponential distribution (the GPD with $\xi = 0$).
  The predictions at various levels of $T$ are shown as the solid line (+) rising to the right. The actual quantiles for the
  sampled distribution are shown by the dashed line and the $1/\sigma$ Bayes predictions (knowing that $\xi = 0$) are shown as circles.
  }\label{extraps}
\end{figure}

A typical extrapolation by the interpolated predictor is illustrated
in Fig.~\ref{extraps}. Extrapolation is from a sample of size $N=15$
drawn from a GPD with $\xi=0$ (i.e. an exponential distribution).  The data
(dots, bottom left) is plotted against the empirical return levels
(both logarithmic), and the predictions are plotted likewise against
the corresponding return level aimed for. The actual
quantiles are also plotted (dashed), together with the $1/\sigma$ Bayes predictions
(knowing that the data is drawn from an exponential).
The predictions of the interpolated predictor typically exceed the actual quantile by a considerable margin, reflecting the uncertainty in the actual
value of $\xi$.
Loosely speaking, although
the data may suggest that there is a high likelihood that the underlying distribution is indeed exponential, the
interpolated predictor recognises that there remains an appreciable chance that the data may have been drawn from
a heavy-tailed GPD with a value of $\xi >0$.

This inherent bias towards larger predictions must not
be confused with risk aversion. The actual
quantiles (dashed lines) can only be drawn here because the
underlying parameters are known. In practice the
parameters will be unknown, and the prediction algorithm
has been designed to allow for the possibility that the data may
have come from any GPD, including those
with heavier tails.

\section{Application to other distributions}
\begin{figure}[hp!] \centering
\includegraphics[width=73mm,keepaspectratio]{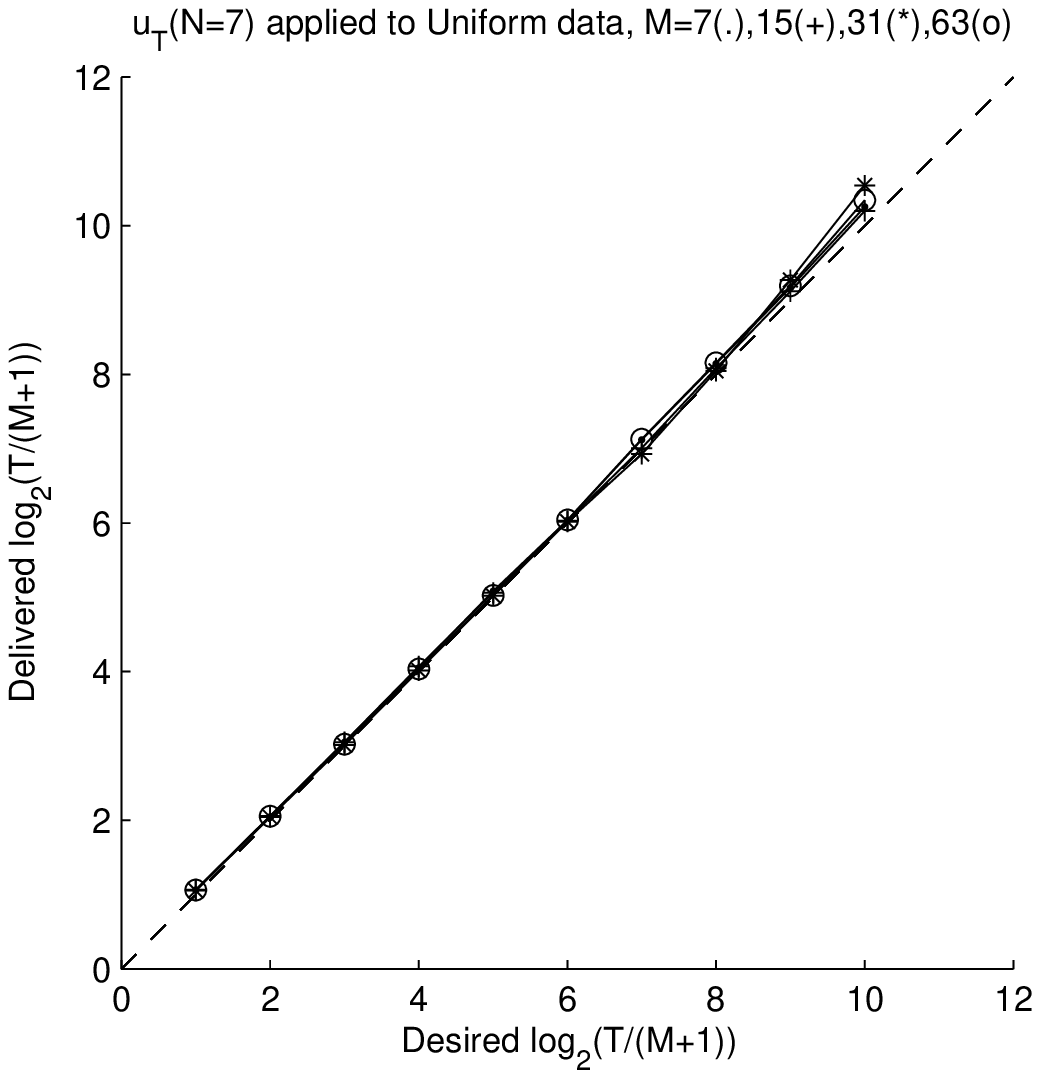}
\includegraphics[width=73mm,keepaspectratio]{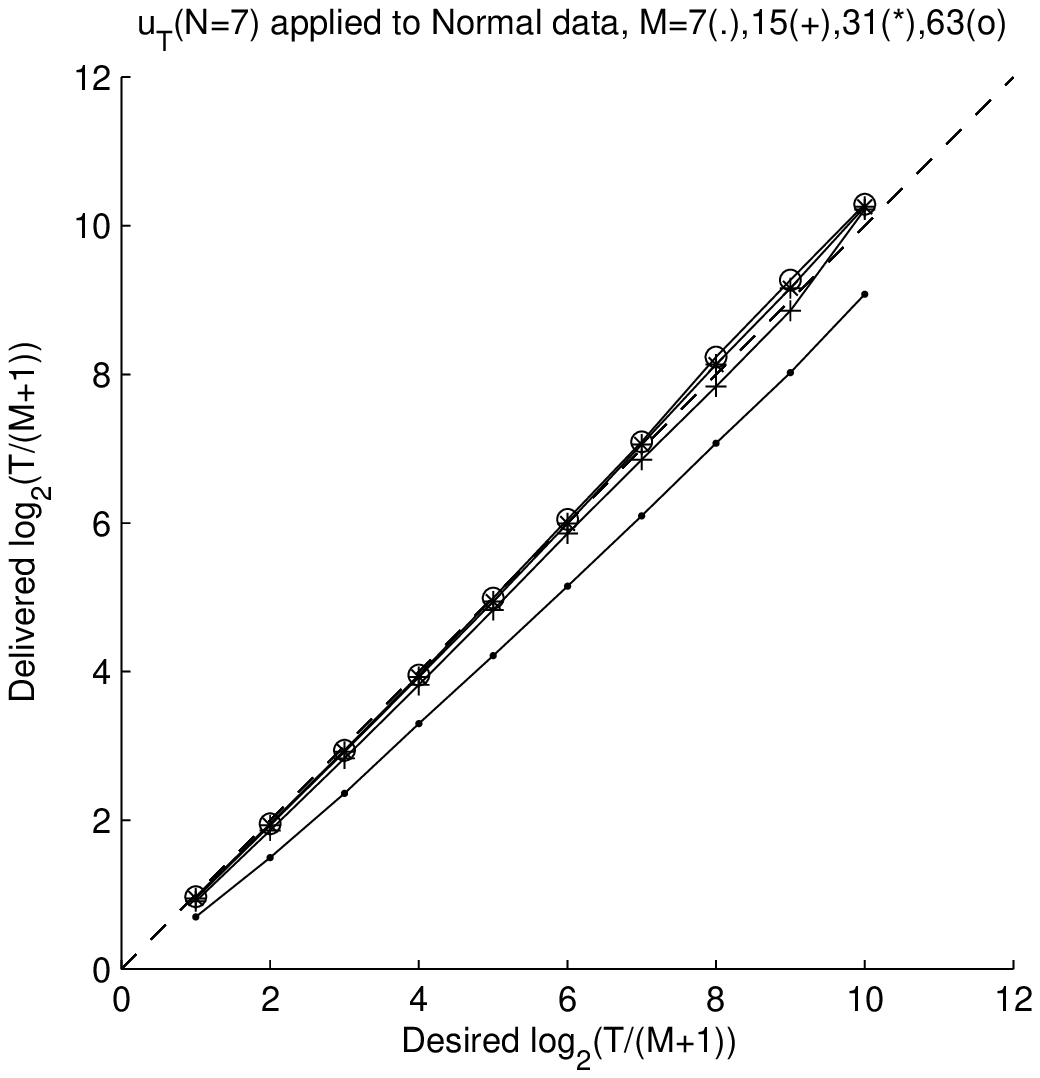}\\
\includegraphics[width=73mm,keepaspectratio]{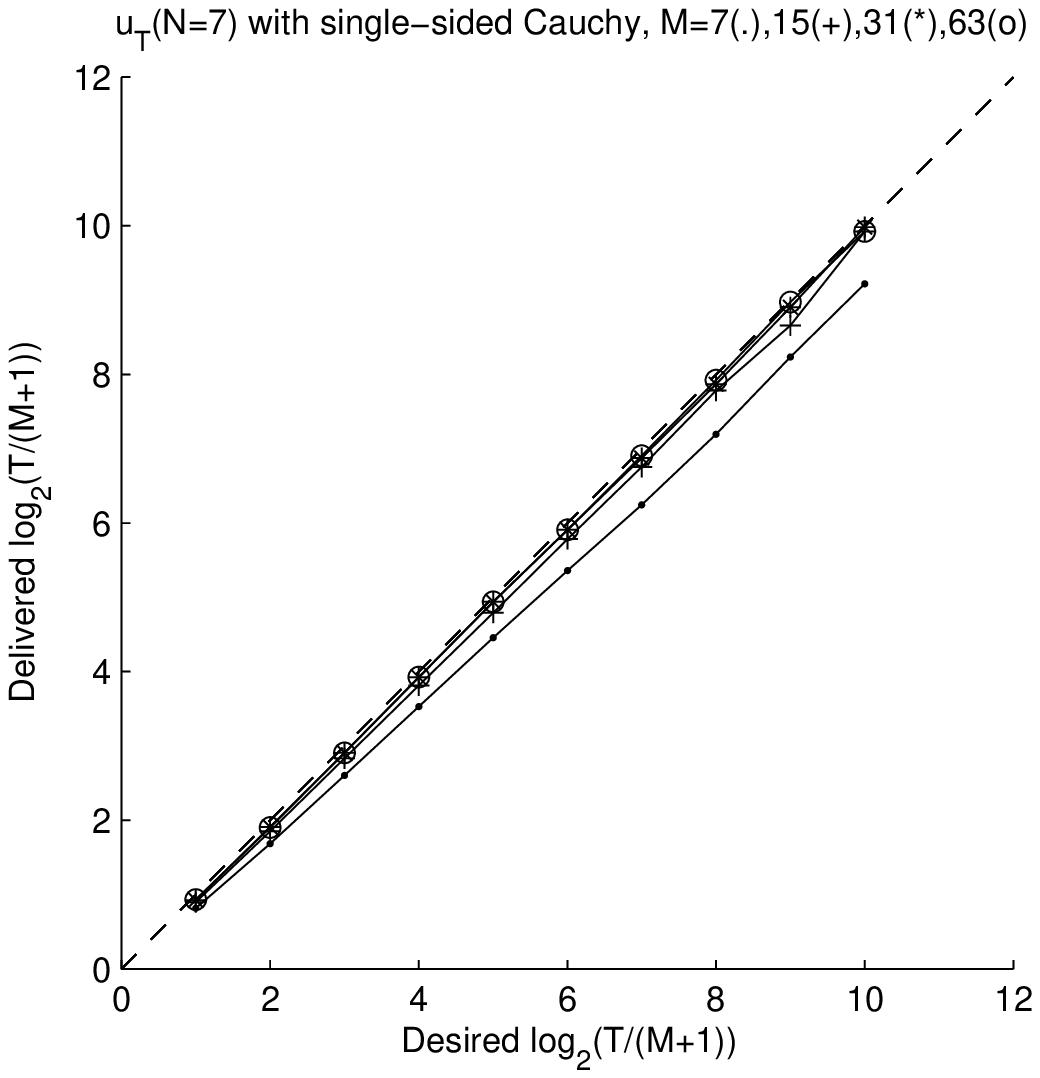}
\includegraphics[width=73mm,keepaspectratio]{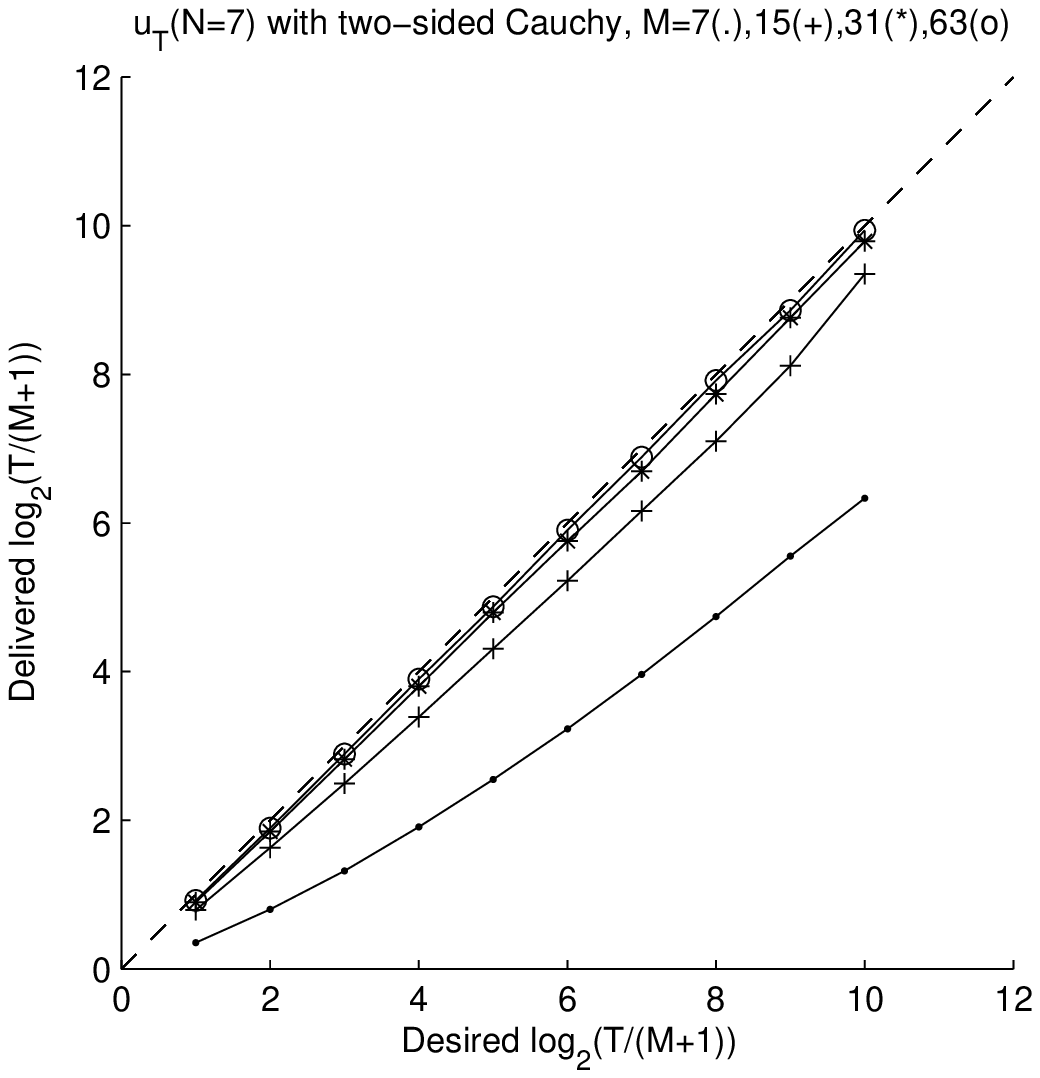}\\
\includegraphics[width=73mm,keepaspectratio]{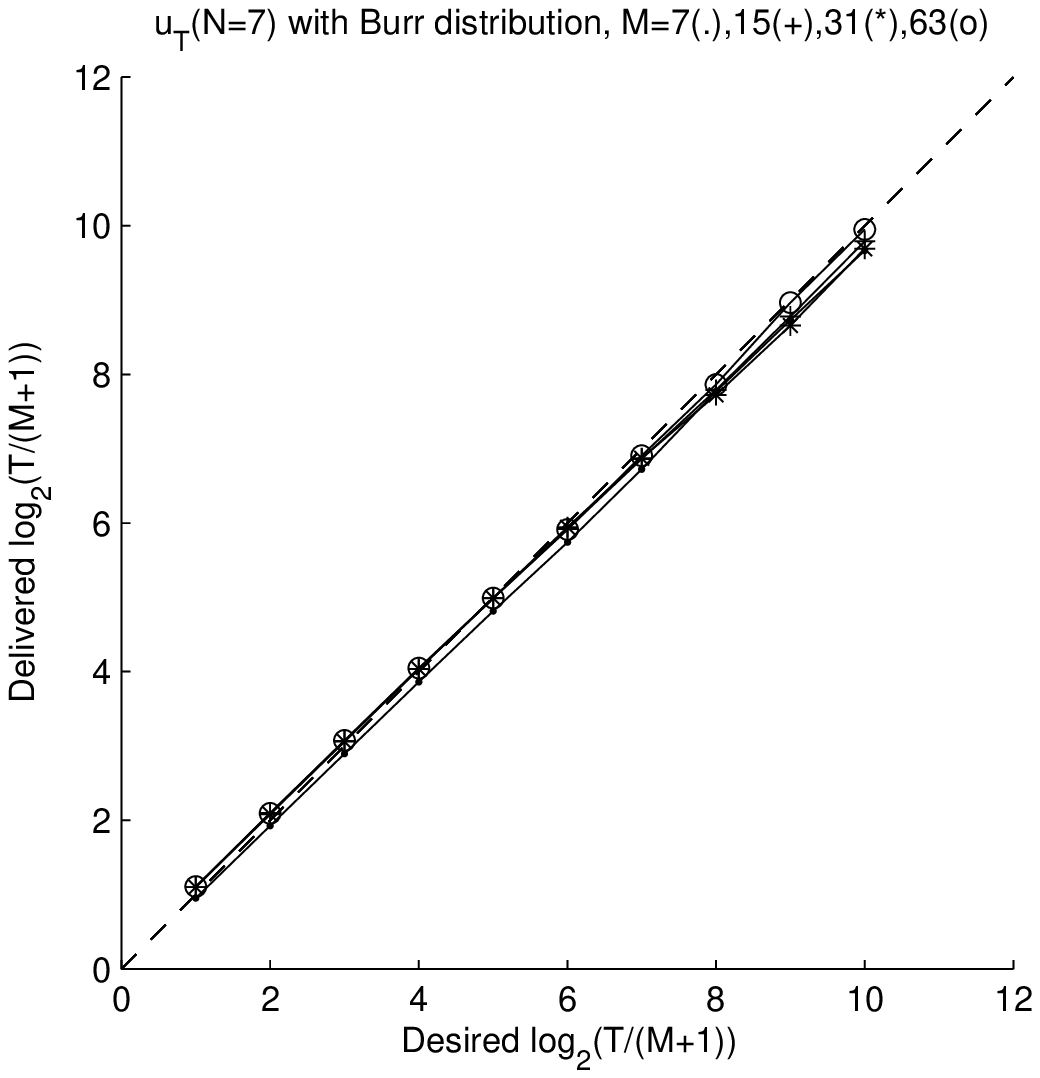}
\includegraphics[width=73mm,keepaspectratio]{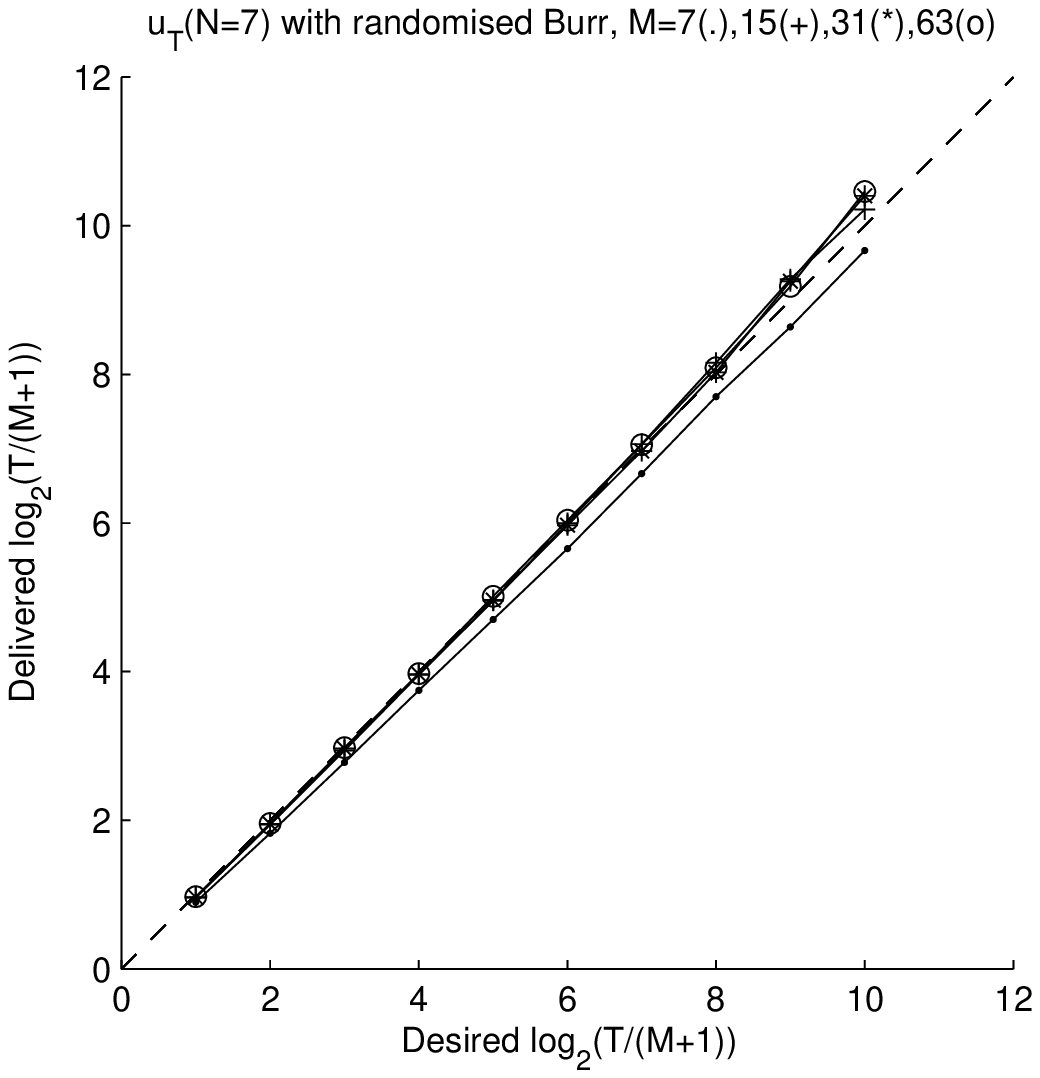}
  \caption{  The probability performance of the interpolated predictor
$u_T$ (Eqn.~\ref{hybrid}) based on the upper $N=7$ data points of
samples of size $M= 7$, 15, 31 and 63 drawn from various
distributions. The distributions are, respectively, uniform, normal,
single-sided Cauchy, two-sided Cauchy, Burr and randomised Burr. In
all cases the location- and scale-parameters are randomised for each
of the $100,000(M+1)$ samples drawn. In the final randomised Burr,
each of the two remaining shape parameters were also randomised.
  }\label{may09trials}
\end{figure}

The interpolated predictor $u_T$ is now applied to samples drawn from a
variety of distributions outside the GPD family. The results are
shown in Fig.~\ref{may09trials}, where the delivered return interval
is plotted against that designed for.

The $N=7$ version of $u_T$ was
used. This is applicable to any sample of size $M \geq 7$, by using
only the upper $N=7$ data points for prediction. The performance is
illustrated for samples drawn from uniform, normal, one- and
two-sided Cauchy and two variants of the Burr distribution.
Both axes of Fig.~\ref{may09trials} plot $\log_2(T/(M+1))$, with $T$ as target on the $x$-axis and as-delivered on the $y$-axis.
The right-most points, with ordinate $\log_2(T/(M+1)) = 10$, thus correspond to extrapolations beyond the data by THREE orders of magnitude - i.e. to the $T \approx 64,000$ level from a sample of size $M = 63$, using just the largest 7 data points thereof.

In all cases, as $M$ becomes significantly larger than $N=7$, the
probability is matched to increasingly better degrees of
approximation. This accords with the general expectation, since the
upper quantiles will approach GPDs as $M \gg N$.

In all cases the location- and scale-parameters were randomised (by
picking $(\mu, \sigma)$ from normal distributions) for each of the
$100,000(M+1)$ samples drawn. However, given that the predictor
guarantees location-scale invariance this should not (and did not)
affect the performance.

The uniform (upper left) delivers good probability matching at
all $M \geq N$, since the uniforms lie within the GPD class.

The normal (upper right) is the first example outside the GPD class.
As may be expected, the predictor based on the $N=7$ upper order statistics
delivers poor performance for a sample of size $M=7$ (because the upper 7 data points are the full data set, and these will have a two-sided normal distribution, and are thus blatantly far from the GPD, and are in no sense extreme values). However, the
probability matching improves dramatically as the sample size $M$ increases even moderately beyond
$N = 7$.

For the two Cauchy examples (centre, Fig.~\ref{may09trials}) the
performance is good in the one-sided case for all $M$, but requires
higher $M$ in the two-sided case. Again, this follows expectations.
It highlights the fact that the problem lies with attempting to
extrapolate from non-extreme data, much of which lies below the mode
in the two-side case. The good results for the one-sided
distribution, even when $M=N$, illustrate that the famously
heavy-tails of the Cauchy present little problem.

The results for the Burr distributions (lower figures), taken from
the family
\begin{equation}
F(y) = 1 - ( 1 + y^A)^B \ \ \ \ A,B >0 \ \ \ \ \mathrm{with} \ y =
\frac{x-\mu}{\sigma}
\end{equation}
again illustrate that the predictor works well with heavy-tailed
distributions. The lower right figure shows the performance over
games played against Burr distributions wherein the parameters $A$
and $B$ (as well as $\mu$ and $\sigma$) are drawn from normal
distributions. The good performance obtained when playing against
such a large stock of distributions with randomised parameters gives
some credence to the claim that $u_T$ has properties approaching
that of a universal probability-matching extreme value predictor.

Good performance has also been found to be delivered for numerous other distributions, 
but it should be noted that there are cases where probabilities do not match well. 
These include Weibull distributions with small shape parameters, 
and beta distributions where the second parameter is small. ``Two population'' distributions
can also be readily constructed for which the predictor performs poorly. These could be said 
to be of the ``black swan'' variety, a simple archetype being 99\% uniform over $[0,a]$ with the remaining 1\% uniform over $[b,1]$, with $a << b$. For moderate $N$, most samples contain no information about the existence of the second population above $b$, the unknown unknowns. Extreme predictions are thus mostly based on data from $[0,a]$ whilst next events with high return period are  values in $[b,1]$. However, for any extrapolation proposal there will obviously be bad cases, and what is surprising about the proposed predictor is how widely and how often it does work (particularly given how small the sample sizes are for which it does work).

\section{Summary}
The intuition that a predictor, designed to match probability within the GPD family, might
transfer some of its potentially-desirable properties across to more general distributions appears to have been borne out.
Although the candidate predictor was only approximately probability-matching across the whole range of GPDs,
the precision with which return levels could be delivered, even from small data sets drawn from non-GPD distributions, is remarkable. 

Finally, it should be emphasised that it is not the intention of this paper to encourage the
extrapolation from three data points to the $T=4096$ level. Rather, suspecting that many current methods of extrapolation
are inherently optimistic, the paper has endeavoured to put forward a novel alternative
for criticism and/or further development.

\section*{Appendix 1: Densities of the normalised data}

Let $N$ raw data points $\mathbf{x} = \lbrace x_1, \ldots, x_N
\rbrace $ be sampled from a GPD with distribution function
\begin{equation}
F(x \ | \ \mu, \sigma, \xi) = 1 - { \left(  1 + \xi \frac{(x-
\mu)}{\sigma } \right) }^{-1 / \xi}
\end{equation}
and density
\begin{equation}
p(x \ | \ \mu, \sigma, \xi) = \frac{1}{\sigma}{ \left(  1 + \xi \frac{(x-
\mu)}{\sigma } \right) }^{-1 -{1 / \xi}}
\end{equation}
These functions exist over the appropriate domains $x>\mu$ for
$\xi>0$ and $\mu <x<\mu-\sigma/\xi$ for $\xi<0$. The case of $\xi =
0$ reduces to the exponential case.

The parameters $(\mu, \sigma, \xi )$ are assumed unknown, and we
consider first the case $\xi>0$, defining $\alpha  = 1/\xi$. This
corresponds to the heavy-tailed case.

We define the ordered data $\mathbf{X} = \mathrm{sort}(\mathbf{x})$,
such that $X_1 \leq X_2 \leq \ldots \leq X_N$.

All possible samples of raw data $\mathbf{x}$ form an N-dimensional
data space. The possible sorted data samples cover only a
semi-infinite prism-shaped subset of a similar N-dimensional space.
The density over that prism of the ordered data space is identical
to that over an equivalent region of the unordered data space
multiplied by a factor $N!$, this being the number of such prisms
required to make up the full space.

The density over the prism-shaped space of ordered data is thus
\begin{equation}
p(\mathbf{X} \ | \ \mu, \sigma,\xi) = N! p(\mathbf{x} \ | \ \mu,\sigma,\xi)
\end{equation}

Although the parameters are unknown, we may define
\begin{equation}
y_j = 1 + \xi \frac{(X_j - \mu)}{\sigma}  \qquad \ \ j = 1, \ldots N
\end{equation}
It follows that $1 \leq y_1 \leq y_2 \ldots \leq y_N $ and over the
domain of  the $y$'s, their density is
\begin{equation}
p(\mathbf{y}|\alpha) = \frac{N! \alpha^{N}}{(y_1 y_2 \ldots y_N
)^{1+\alpha}}
\end{equation}

We define $(N-2)$ location- and scale-independent statistics
$\mathbf{t}$ via
\begin{equation}
t_j \equiv \frac{X_{j+1} - X_1}{X_N - X_1} \qquad j=1,\ldots (N-2)
\end{equation}
This is the normalised data $\mathbf{t}$. It is ordered and lies in
the unit interval ( $0 \leq t_1 \leq t_2 \leq \ldots \leq t_{N-2}
\leq 1$). Trivially, we can also measure from the opposite end of
the interval, defining $\tau_j = 1 - t_j$

We also define
\begin{equation}
q = \frac{y_N - y_1}{y_N}
\end{equation}

Writing the left and right end points of the unit interval as $t_0 =
0$, $\tau_0 =1$   and  $t_{N-1} = 1$, $\tau_{N-1} =  0$ respectively
then together the $N$ transformations
\begin{equation}
y_j = y_N (1 - q \tau_j)   \ \ \ \mathrm{for} \ j = 0, \ldots , (N-1)
\end{equation}
have Jacobian $y_1^{N-1} q^{N-2}$ and the density in the new
variables is
\begin{equation}
p(y_N,q,\btau|\alpha) = \frac{N! \alpha^{N-1} y_N^{-N\alpha -1}
q^{N-2} }{\prod_{j=0}^{N-2} (1-\tau_j q )^{1+\alpha} }
\end{equation}

For $q$ fixed, the variable $y_N$ is bounded below at $1/(1-q)$. It
is removed by integration, using
\begin{equation}
\int_{(1-q)^{-1}}^{\infty}y_N^{-N\alpha -1} \ dy_N =
\frac{(1-q)^{N\alpha}}{N\alpha}
\end{equation}
before removal of $q$ via the integration
\begin{equation}
p(\bt|\alpha) = (N-1)! \ \alpha^{N-1} \int_0^1 q^{N-2}
(1-q)^{(N-1)\alpha -1} \prod_{j=1}^{N-2} (1- \tau_j q)^{-(1+\alpha)}
\ dq
\end{equation}

This is a standard Euler integral representation of a Lauricella
function $F_D^{(N-2)}$ \citep{extonred} 
leading to the desired density
\begin{equation}
p(\bt|\alpha) = N_1! \ \alpha^{N_1}
 \ \Gamma \left[
\begin{array}{c} N_1, \ N_1 \alpha \\ N_1(1+\alpha) \end{array} \right]
F_D^{(N-2)}\left( N_1, \bone +\balpha; N_1(1+\alpha); \boldtau
\right) \label{poft1}
\end{equation}
where $N_1 = N-1$, and $\bone + \balpha$ is a vector of length $N-2$
with each element being $1+\alpha$.

The Lauricella function may be expressed using one of its Euler
transformations
\citep{extonblue2} 
to give the alternative expression
\begin{equation}
p(\bt|\alpha) = N_1! \ \alpha^{N_1}
 \ \Gamma \left[
\begin{array}{c} N_1, \ N_1 \alpha \\ N_1(1+\alpha) \end{array} \right]
\frac{t_1^{(N-2)\alpha -1}}{(t_2 \ldots t_{N-2})^{1+\alpha}}
F_D^{(N-2)}\left( N_1 \alpha, \bone +\balpha; N_1(1+\alpha);
\boldkappa \right) \label{alphadens}
\end{equation}
where
\begin{equation}
\kappa_1 = 1- t_1, \ \ \qquad \kappa_j = 1 - \frac{t_1}{t_j} \ \ \
\mathrm{for} \ j = 2, \ldots N-2
\end{equation}
The advantage of this representation is that for sufficiently small
$\alpha$, the leading Lauricella parameter $(N-1)\alpha$ is small
such that the Lauricella function is close to unity.

For the bounded-tail case, we define $\beta \equiv -1/\xi$ with
$\beta>0$, and a similar derivation leads to
\begin{equation}
p(\bt|\beta) = N_1! \ \beta^{N_1}
 \ \Gamma \left[
\begin{array}{c} N_1, \ \beta \\ N_1 + \beta \end{array} \right]
F_D^{(N-2)}\left( N_1, \bone -\bbeta; N_1+\beta; \bt \right)
\label{betadens}
\end{equation}
where $\bbeta = \beta \bone$. It should be noted that this is not
simply the $\xi >0$ expression with the substitution $\alpha
\rightarrow -\beta$ since the bounded nature of the tails leads 
to different limits in the various integrations.

\section*{Appendix 2. Derivation of $G(s_T(\bt)|\xi)p(\bt|\xi)$ }
The integral
\begin{equation}
\int_{\forall \bt} \ G(s_T(\bt)|\xi) \ p(\bt|\xi) \ d\bt
\label{app2eq1}
\end{equation}
is the ({\em average}) probability performance delivered at some
given tail parameter $\xi$ by any (normalised) predictor $s_T(\bt)$
which is a function of the (normalised) data $\bt$. The rest of this
paper concerns itself with attempting to construct such a function
$s_T(\bt)$ such that the integral \ref{app2eq1} is independent of
$\xi$.

The functional form of the integrand $G(s_T(\bt)|\xi)p(\bt|\xi)$ can
be derived from first principles in a manner akin to the derivations
of Appendix 1 or, equivalently, via appropriate integration of the
Appendix 1 results, as here.


Consider an ordered sample of size $(N+1)$ drawn from a GPD with
shape parameter $\xi$. Let the normalised data be $\bh = \lbrace
h_1, \ldots h_{N-1} \rbrace$. This normalised data is related to the
normalised data $\bt$ plus an extra normalised data point $s$ via
$h_j = t_j/s$ for $j=1$ to $N-1$.

For $\xi > 0$  ($\xi = 1/\alpha$), the density of $\bh$ is
obtained from Eqn.~(\ref{poft1}) as
\begin{equation}
p(\bh|\alpha) = \frac{N!}{N+1} \ \alpha^{N}
 \ \Gamma \left[
\begin{array}{c} N, \ N \alpha \\ N(1+\alpha) \end{array} \right]
F_D^{(N-1)}\left( N, \bone +\balpha; N(1+\alpha); \bone-\bh \right)
\label{pofh1}
\end{equation}
The $N+1$ divisor must be introduced since we are considering only
that $1/(N+1)$ fraction of cases where the next data point $x_{N+1}$
exceeds the historical data maximum.

Substituting $h_j = t_j/s$ (a transformation with Jacobian $1/s^N$)
leads to the joint density for $\bt$ and $s>1$ of
\begin{equation}
p(s,\bt|\alpha) = \frac{N!}{N+1} \ \alpha^{N}
 \ \Gamma \left[
\begin{array}{c} N, \ N \alpha \\ N(1+\alpha) \end{array} \right]
\frac{1}{s^N} F_D^{(N-1)}\left( N, \bone +\balpha; N(1+\alpha);
\bone-\frac{\bt}{s} \right) \label{pofh1}
\end{equation}
where the vectors extend from $j=1$ to $N-1$.

Using one of the standard Lauricella transforms \citep{extonblue2}
this can be written 
\begin{equation}
p(s,\bt|\alpha) = \frac{N!}{N+1} \ \alpha^{N}
 \ \Gamma \left[
\begin{array}{c} N, \ N \alpha \\ N(1+\alpha) \end{array} \right]
F_D^{(N-1)}\left( N, \bone +\balpha; N(1+\alpha); 1-s, 1-\bt \right)
\end{equation}
This now needs to be integrated from $s = s_T$ to $\infty$ to obtain
the tail probability.

It follows from the Euler integral representation of the Lauricella
function that
\begin{eqnarray}
\int_{s_T}^{\infty} \
F_D^{(n)}(a,b_1,\ldots,b_n;c;1-s,x_2,\ldots,x_n) \ ds \hspace*{40mm} & & \nonumber \\
=
\frac{(c-1)}{(b_1-1)(a-1)}F_D^{(n)}(a-1,b_1-1,b_2,\ldots,b_n;c-1;1-s_T,x_2,\ldots,x_n)
& &
\end{eqnarray}

 This
leads to
\begin{eqnarray}
G(s_T(\bt)|\alpha)p(\bt|\alpha) & = & \int_{s_T}^{\infty} \ p(s,\bt |\alpha) \ ds \\
 & =  & \frac{N!}{N+1}  \ \alpha^{N-1}
 \ \Gamma \left[
\begin{array}{c} N-1, \ N \alpha \\N+N\alpha -1  \end{array}
\right] \ldots \nonumber \\
& & \times \ F_D^{(N-1)}\left( N-1, \alpha , \bone +\balpha ;
N+N\alpha-1; 1-s, 1-\bt \right)
\end{eqnarray}

A further Lauricella transform leads to the final form
\begin{eqnarray}
G(s_T(\bt)|\alpha)p(\bt|\alpha)  & = & \frac{N!}{N+1} \ \alpha^{N-1}
 \ \Gamma \left[
\begin{array}{c} N-1, \ N \alpha \\ N+N\alpha -1 \end{array} \right]
\frac{t_1^{(N-1)\alpha -1}}{(t_2 \ldots
t_{N-2})^{1+\alpha}}\frac{1}{s^\alpha} \ldots  \nonumber \\
& & \times F_D^{(N-1)}\left( N\alpha, \bone +\balpha , \alpha;
N+N\alpha-1; \br \right)\label{gstptalf}
\end{eqnarray}
with $r_j = 1-t_1/t_j$ for $j=2$ to $N$ (where $t_{N-1} =1$ and $t_N
= s_T$) .

A similar derivation, omitted for brevity, for the bounded-tail case
with $\xi$ negative ($\beta = -1/\xi$ positive), leads to
\begin{eqnarray}
G(s_T(\bt)|\beta)p(\bt|\beta)  & = & \frac{N!}{N+1} \ \beta^{N_1}
 \ \Gamma \left[
\begin{array}{c} N_1, \ 1+ \beta \\ N+\beta \end{array} \right]
\ldots \nonumber \\ & & \times \frac{1}{s_T^{N_1}}F_D^{(N_1)}\left(
N_1, \bone - \bbeta; N+\beta; \frac{\bt}{s_T} \right)
\label{gstptbet}
\end{eqnarray}
where $\bt = \lbrace t_j \rbrace$ for $j = 1$ to $N_1=N-1$ (with
$t_{N-1} = 1$). Again, this is not simply the $\xi >0$  expression
with $\alpha \rightarrow -\beta$, owing to the bounded nature of the
tails changing various integration limits.

Eqns.~\ref{gstptalf} and \ref{gstptbet} for the exceedance
probability of the predictor play a central role in the main paper.
For $s_T(\bt)$ to be a probability-matching predictor, the
integral of $G(s_T(\bt)|\xi)p(\bt|\xi)$ over all possible data $\bt$
should equal the desired exceedance probability $1/T$, and this
should be true at any value of the tail parameter $\xi$.

\section*{Appendix 3: The constraint on the power law exponents in the bounded-tail case}

We derive forms for power law predictors for the extreme bounded-tail
case (where $\xi = -1/\beta$ and $\beta$ is small and positive).

In the $\xi > 0$  case, a constraint on the power-law exponents was
derived directly from the small $\alpha$ behaviour of the Lauricella
form of $G(s_T(\bt)|\alpha)p(\bt|\alpha)$ of Eqn.~\ref{eqnpp1}.
Unfortunately, a similar procedure does not seem possible via
$G(s_T(\bt)|\beta)p(\bt|\beta)$ of Eqn.~\ref{eqnpp2} (or any of its
Lauricella transformations or asymptotic forms) in the bounded-tail case.
Instead, it is necessary to return to first principles, inserting
the small $\beta$ approximation en route.

For a sample of size $N$ drawn from a GPD in the bounded-tail case, the
density of the ordered, normalised data $\bt$ is given by a
Lauricella transformation of Eqn.~\ref{Frechdens} as
\begin{eqnarray}
p(\bt |\beta) \ d\bt & =&  N_1! \ \beta^{N_1} \ \Gamma \left[
\begin{array}{c} N_1, \ \beta \\ N_1+\beta \end{array} \right]
\frac{\tau_{N-2}^{2\beta-1}}{(\tau_1 \ldots \tau_{N-3})^{1-\beta}}
\nonumber  \\ & &  \hspace*{10mm} \times \ F_D^{(N-2)}(\beta,
N_1\beta+1, 1-\bm{\beta};N_1+\beta;\bm{\Psi}) \ d\bt \label{app3eq1}
\end{eqnarray}
with $N_1 = N-1$, $\bm{\tau} = 1 - \bm{t}$, and $\psi_j = 1 -
\tau_{N-2}/\tau_j$ for $j = 0$ to $N-3$ (where $\tau_0 = 1$).

Consider now drawing an additional data point $x_{N+1}$. This might
not exceed the sample maximum, but we restrict attention to the
$1/(N+1)$ fraction of cases when it does. The density of the
ordered, normalised data $\bm{t^\dagger}$ in this region is thus
given by Eqn.~\ref{app3eq1}, with $N \rightarrow N+1$ and $\bm{t}
\rightarrow \bm{t}^\dagger = [\bm{t}/s_*, 1/s_*]$. That is,
\begin{eqnarray}
p( {\bt}^{\dagger} | \beta ) \ d\bt^{\dagger} & = &  \frac{N!}{N+1}
 \ \beta^{N} \ \Gamma \left[
\begin{array}{c} N, \ \beta \\ N+\beta \end{array} \right]
\frac{(\tau_{N-1}^{\dagger} )^{2\beta-1}}{(\tau_1^{\dagger} \ldots
\tau_{N-2}^{\dagger} )^{1-\beta}} \nonumber  \\
 & &   \hspace*{10mm} \times \ F_D^{(N-1)}(\beta,
N\beta+1, 1-\bm{\beta};N+\beta;\bm{\Psi}^{\dagger} ) \
d\bt^{\dagger}  \label{app3eq2}
\end{eqnarray}
where
\begin{equation}
\tau_j^\dagger = 1- t_j^\dagger = 1 - t_j/s_* = \frac{u_*+\tau_j}{1+
u_*} \ \ \mathrm{and} \ \psi_j^\dagger = 1-
\frac{\tau_{N-1}^\dagger}{\tau_j^\dagger} = \frac{\tau_j}{u_* +
\tau_j}
\end{equation}
with $u_* = s_* -1$. The leading $N+1$ in the denominator of
Eqn~\ref{app3eq2} accounts for the fact we only consider those cases
where the next data point exceeds the sample maximum.

The Jacobian of the transformation gives
\begin{equation}
d\bt^\dagger = \frac{1}{s_*^N}d\bt ds_* =
\frac{1}{(1+u_*)^N}d\bm{\tau}du_*
\end{equation}

leading to
\begin{eqnarray}
& & p( \btau, u_* | \beta) \  d\btau du_* \ = \  \frac{N!}{N+1}  \
\beta^{N} \ \Gamma \left[
\begin{array}{c} N, \ \beta \\ N+\beta \end{array} \right]
\frac{u_*^{2\beta-1}}{(1+u_*)^{N\beta+ 1} \prod_{j =
1}^{N-2}(u_*+\tau_j)^{1-\beta}} \nonumber  \\
 & &   \hspace*{10mm} \times \ F_D^{(N-1)}(\beta,
N\beta+1, 1-\bm{\beta};N+\beta; \frac{1}{1+u_*},\frac{\btau}{u_* +
\btau} ) \ d\btau du_* \label{app3eq3}
\end{eqnarray}

We know that
\begin{equation}
d\btau \int_{u(\btau)}^\infty \ p(\btau, u_* | \beta) du_* =
G(u(\btau) |\beta)p(\btau | \beta) d\btau \
\end{equation}
and already have an exact expression for this (Eqn \ref{gstptbet}).
However, it has not yet been possible to deduce a predictor, not
even an approximate one, from that expression. Instead we need to
pursue a more circuitous route.

Rather than considering the tail integral directly (i.e. $G(u
(\btau)|\xi)p(\btau|\xi)$), we endeavour to approximate instead the
integral over the region between the data maximum and the
prediction, viz:
\begin{equation}
\int_{\forall \btau} \ \int_0^u \ p(\btau, u_*) \ du_* \ d\btau =
\frac{1}{N+1} - \frac{1}{T}
\end{equation}
leading to
\begin{eqnarray}
& & 1-\frac{N+1}{T} = N! \beta^{N}  \Gamma \left[
\begin{array}{c} N, \ \beta \\ N+\beta \end{array} \right]
\int_{\forall \btau} \ \int_0^u
\frac{u_*^{2\beta-1}}{(1+u_*)^{N\beta+ 1} \prod_{j =
1}^{N-2}(u_*+\tau_j)^{1-\beta}} \nonumber  \\
 & &   \hspace*{10mm} \times \ F_D^{(N-1)}(\beta,
N\beta+1, 1-\bm{\beta};N+\beta; \frac{1}{1+u_*},\frac{\btau}{u_* +
\btau} ) \ du_* d\btau  \label{app3eq4}
\end{eqnarray}
This expression is exact, and if a function $u(\tau)$ could be found
which satisfies this then we would have our exact probability
matching predictor for all $\beta>0$. This has not yet been
possible, and thus we proceed with an approximation.

At the first level of approximation, we could set
$F_D^{(N-1)}(\beta, \ldots) \approx 1$ and for $u \ll \tau_j$ we
could ignore the $u_*$ terms in the denominator. Postulating a power
law predictor of the form $u = \prod_j \tau_j^{\rho_j}$ then leads
to a constraint equation for the exponents $\rho_j$. However, the
resulting predictor has very limited validity: it applies at $N=3$,
but for higher $N$ it is valid only for very small extrapolations
beyond the data. We thus proceed directly to an improved
approximation.

Assuming $F_D^{(N-1)}(\beta, \ldots) \approx 1$, and writing $u_* =
w u$ the inner integral becomes (for $u \ll 1 $)
\begin{eqnarray}
I_u & = & \int_0^u \ \frac{u_*^{2\beta-1}}{(1+u_*)^{N\beta+ 1}
\prod_{j = 1}^{N-2}(u_*+\tau_j)^{1-\beta}} \ du_* \\
    & = & u^{2\beta-1}  \prod_{j=1}^{N-2} \tau_j^{\beta -1} . \int_0^u
    w^{2\beta -1} \left[ \prod_{j=1}^{N-2} \left( 1 + \frac{u}{\tau_j}w
    \right)^{\beta -1} \right] \ (1+uw)^{-N\beta-1} \ u dw \nonumber \\
    & = & u^{2\beta}  \prod_{j=1}^{N-2} \tau_j^{\beta -1} \Gamma \left[
\begin{array}{c} 2\beta, \ 1 \\ 1+2\beta \end{array} \right]
F_D^{(N-1)}(2\beta, 1+N\beta, 1-\bm{\beta};1+2\beta;-u,
\frac{-u}{\btau}) \nonumber
\end{eqnarray}

For small $\beta$, the gamma function approaches $1/2\beta$.

The Lauricella function can be transformed on its argument of
largest magnitude
 $-u/\tau_{N-2} \equiv x_1$ using \citet{extonblue2} (Eqn.
4.2.4, second transformation). The new arguments are
\begin{equation}
\frac{x_1}{x_1-1} = \frac{u}{u+\tau_{N-2}} \ \ \mathrm{and} \ \
\frac{x_1 - x_j}{x_1-1} = \left( \frac{u}{u+\tau_{N-2}} \right)
\left( 1-\frac{\tau_{N-2}}{\tau_j} \right)
\end{equation}
These are of interest only in that they are less than unity, and
thus the transformed Lauricella function $F_{D}^{(N-1)}(2\beta,
\ldots)$ can be approximated as unity for $\beta$ small. This only
leaves the leading factor in the Lauricella transformation
\begin{equation}
(1-x_1)^{-a} = \left( 1+\frac{u}{\tau_{N-2}}\right)^{-2\beta}
\end{equation}
and we thus obtain
\begin{equation}
I_u \approx \frac{1}{2\beta}\left(
\frac{u}{u+\tau_{N-2}}\right)^{2\beta} \ \tau_{N-2}^{2\beta} \
\prod_{j=1}^{N-2}\tau_j^{\beta -1}
\end{equation}
Now, assume a power-law predictor such that
\begin{equation}
\frac{u}{u+\tau_{N-2}} = \prod_{j=1}^{N-2}\tau_j^{\rho_j} \ \
\mathrm{for} \ \mathrm{some} \ \rho_j
\end{equation}

To match probability in the extreme bounded-tail limit $\xi \rightarrow -\infty$ we require
\begin{equation}
1-\frac{N+1}{T} = \frac{N!}{2} \ \beta^{N-2} \int_{\forall \btau} \
\tau_{N-2} \prod_{j=1}^{N-2} \tau_j^{\beta(2\rho_j+1)-1} \ d\btau
\end{equation}
(where we have approximated the two gamma functions by $1/\beta$ and
$1/2\beta$ as appropriate for $\beta$ small).

The domain of integration is $\int_0^1 d\tau_1 \int_0^{\tau_1}
d\tau_2 \ldots \int_0^{\tau_{N-3}} d\tau_{N-2}$, leading to

\begin{eqnarray}
1-\frac{N+1}{T}  & = &  \frac{N! \beta^{N-2}}{2. (2 \rho_{N-2}+3)
\beta
\ldots (2 ( \rho_{N-2}+\ldots+\rho_1 ) + N )\beta } \\
&   =  & \frac{N!}{2} \prod_{k=1}^{N-2} \left( 2\eta_k +
(k+2)\right)^{-1}
\end{eqnarray}

That is, our predictor $u(\btau)$ is
\begin{equation}
u = \tau_{N-2} \ \frac{\tauprod}{1-\tauprod} \equiv u_\beta
\end{equation}
where the exponents satisfy
\begin{equation}
\prod_{k=1}^{N-2} \left[ \frac{2\eta_k}{k+2}+1 \right] =
\frac{1}{1-\frac{N+1}{T}} \label{Frech1}
\end{equation}
with partial sums
\begin{equation}
\eta_k = \sum_{i = N-1-k}^{N-2} \rho_i \label{Frech2}
\end{equation}
The above derivation suggests that this should give the correct
exceedance probability $1/T$ in the extreme bound-tail limit of GPDs
with $\xi$ large and negative, and the numerical results presented
in the main text (Fig.~\ref{threebounded}) suggest that it does.

\bibliography{gpdbib}

\begin{thebibliography}{9}
\expandafter\ifx\csname natexlab\endcsname\relax\def\natexlab#1{#1}\fi
\expandafter\ifx\csname url\endcsname\relax
  \def\url#1{\texttt{#1}}\fi
\expandafter\ifx\csname urlprefix\endcsname\relax\def\urlprefix{URL }\fi

\bibitem[{Coles and Tawn(2005)}]{colestawn}
Coles, S., Tawn, J., 2005. Bayesian modelling of extreme surges on the {UK}
  east coast. Phil. Trans. R. Soc. A 363, 1387--1406.

\bibitem[{Datta and Mukerjee(2004)}]{datta1}
Datta, G.~S., Mukerjee, R., 2004. Probability Matching Priors: Higher Order
  Asymptotics. Lecture Notes in Statistics. Springer, New York.

\bibitem[{Embrechts et~al.(1999)Embrechts, Kl\"{u}ppelberg, and
  Mikosch}]{Embrechts}
Embrechts, P., Kl\"{u}ppelberg, C., Mikosch, T., 1999. Modelling Extreme Events
  for Insurance and Finance. Springer, Berlin.

\bibitem[{Exton(1976)}]{extonblue2}
Exton, H., 1976. Multiple Hypergeometric Functions and Applications. Ellis
  Horwood, Chichester, UK.

\bibitem[{Exton(1978)}]{extonred}
Exton, H., 1978. Handbook of Hypergeometric Integrals. Ellis Horwood,
  Chichester, UK.

\bibitem[{McRobie(2004)}]{mcrobieEEE}
McRobie, F.~A., 2004. Exact exceedance estimators and the $1/\sigma$ reference
  prior. Tech. Rep. CUED/D-STRUCT/TR213, Cambridge University Engineering Dept.

\bibitem[{McRobie(2013{\natexlab{a}})}]{McRobieGEV}
McRobie, F.~A., 2013{\natexlab{a}}. Elemental estimators for the {G}eneralized
  {E}xtreme {V}alue tail. arxiv:1304.4362.

\bibitem[{McRobie(2013{\natexlab{b}})}]{McRobieGPD}
McRobie, F.~A., 2013{\natexlab{b}}. Elemental unbiased estimators for the
  {G}eneralized {P}areto tail. arxiv:1304.3918.

\bibitem[{Sweeting(2008)}]{sweeting1}
Sweeting, T., 2008. On predictive probability matching priors. In: Clarke, B.,
  Ghosal, S. (Eds.), IMS Collections: {P}ushing the {L}imits of {C}ontemporary
  {S}tatistics: {C}ontributions in {H}onor of {J}ayanta {K}. {G}hosh. No.~3.
  pp. 46--59.

\end{thebibliography}

\end{document}